\documentclass[12pt]{amsart}
\usepackage{amssymb, amscd, pb-diagram}
\usepackage{verbatim}
\usepackage{hhline}
\usepackage{xspace}

\theoremstyle{plain}
\newtheorem{thm}{Theorem}[section]
\newtheorem{prop}[thm]{Proposition}
\newtheorem{lma}[thm]{Lemma}
\newtheorem{lemma}[thm]{Lemma}
\newtheorem{cor}[thm]{Corollary}

\theoremstyle{definition}
\newtheorem{definition}[thm]{Definition}
\newtheorem{defn}[thm]{Definition}

\theoremstyle{remark}
\newtheorem*{remark}{Remark}
\newtheorem*{remarks}{Remarks}
\newtheorem{example}{Example}

\newcommand{\Inv}{\operatorname{Inv}}
\newcommand{\rank}{\operatorname{rank}}

\newcommand{\Pisom}{\operatorname{Pisom}}
\newcommand{\alglim}{\underrightarrow{\operatorname{alg\,lim}}\,}
\newcommand{\row}{\operatorname{row}}
\newcommand{\col}{\operatorname{col}}
\newcommand{\ds}{\displaystyle}
\newcommand{\tph}{\tilde{\phi}}
\newcommand{\tps}{\tilde{\psi}_1}
\newcommand{\AAA}{\mathcal{A}}
\newcommand{\bbC}{\mathbb{C}}

\newcommand{\bbz}{\mathbb{Z}}
\newcommand{\GG}{\mathcal{G}}

\newcommand{\sF}{{\mathcal F}}
\newcommand{\sA}{{\mathcal A}}
\newcommand{\cstar}{\ensuremath{\text{C}^{*}}\nobreakdash-\hspace{0 pt}}
\renewcommand{\star}{\ensuremath{{}^{*}}\nobreakdash-\hspace{0 pt}}
\newcommand{\cal}{\mathcal}

\def\Sys{{\mathrm{Sys}\thinspace}}

\newenvironment{pf}
   {\begin{proof}}
   {\end{proof}}

\newcommand{\oc}{\text{oc}}  
\newcommand{\op}{\text{op}}  
\newcommand{\loc}{\text{loc}} 
\newcommand{\lop}{\text{lop}}  
\newcommand{\reg}{\text{reg}}
\newcommand{\orcon}{order conserving\xspace}  
\newcommand{\orpre}{order preserving\xspace}  
\newcommand{\lorcon}{locally order conserving\xspace}
\newcommand{\lorpre}{locally order preserving\xspace}
\newcommand{\Orcon}{Order conserving\xspace} 
\newcommand{\Orpre}{Order preserving\xspace}

\newcommand{\orconion}{order conservation\xspace}
\newcommand{\lorconion}{local order conservation\xspace}
\newcommand{\orpreion}{order preservation\xspace}
\newcommand{\lorpreion}{local order preservation\xspace}

\newcommand{\Orconion}{Order conservation\xspace}

\newcommand{\vstrl}{\rule[-2pt]{0 pt}{16.5 pt}}  
\newcommand{\vstrs}{\rule[-2pt]{0 pt}{15 pt}}  

\newcommand{\nofusion}[4]{\ensuremath{
\begin{array} {||c||c||}
\hhline{|t:=:t:=:t|} 
\vstrs #1 & #4 \\
\hhline{|:=::=:|} 
\vstrs #3 & #2 \\
\hhline{|b:=:b:=:b|} 
\end{array} 
}}

\newcommand{\hfusion}[4]{\ensuremath{
\begin{array} {||cc||}
\hhline{|t:==:t|}
\vstrs #1 & #4 \\
\hhline{|:==:|}
\vstrs #3 & #2 \\
\hhline{|b:==:b|}
\end{array}
}}

\newcommand{\vfusion}[4]{\ensuremath{
\begin{array} {||c||c||}
\hhline{|t:=:t:=:t|}
\vstrl #1 & #4 \\
\vstrl #3 & #2 \\
\hhline{|b:=:b:=:b|}
\end{array}
}}

\newcommand{\fullfusion}[4]{\ensuremath{
\begin{array} {||cc||} 
\hhline{|t:==:t|}
\vstrl #1 & #4 \\
\vstrl #3 & #2 \\
\hhline{|b:==:b|}
\end{array}
}}

\numberwithin{equation}{section}
\begin{document}
\title[Limits of nest algebras]{Limits of 
finite dimensional nest algebras}
\date{November 4, 2000}
\author{Alan Hopenwasser}
\address{Dept. of Mathematics \\  University of Alabama \\
 Tuscaloosa, AL 35487 \\ U.S.A.}
\email[Alan Hopenwasser]{ahopenwa@euler.math.ua.edu} 
\author{Stephen~C. Power}
\address{Dept. of Mathematics \& Statistics \\ Lancaster
 University \\ Lancaster, U.K. LA1 4YF}
\email[Stephen~C. Power]{s.power@lancaster.ac.uk}
\thanks{The first author was supported by an EPSRC visiting fellowship
 to Lancaster University}
\thanks{The second author was supported by an EPSCoR of Alabama 
 travel grant to visit the University of Alabama.}

\begin{abstract}
We introduce order conserving embeddings as a more general form of
order preserving embeddings
between finite dimensional nest algebras.
The structure of these embeddings is determined, in terms of order
indecomposable decompositions, and they are shown to be determined up to
inner conjugacy by their induced maps on $K_0$.
Classifications of direct systems and limit algebras are obtained in
terms of dimension distribution groups.
\end{abstract}
\maketitle

\tableofcontents

\newpage

\section{Introduction} \label{S:intro}

One of the major themes in the study of direct systems of operator
algebras is the classification of such systems and their limit
algebras.  Since universal results are largely out of reach and
inappropriate to specific families,
 the usual approach is to study families of systems which
are simultaneously tractable and of intrinsic interest and, for each
such family, to find a complete set of invariants (up to a suitable
notion of isomorphism).  The selection of a family for classification
entails the specification of a family of (usually finite dimensional)
``building block'' algebras and a class of allowable embeddings
between the building block algebras.  One desirable goal in the
overall process is to achieve greater unification of the various
classifying invariants.

The first example of this scheme was Glimm's classification of UHF
\cstar algebras (building blocks: full matrix algebras 
($M_n(\bbC)$); embeddings: unital \star homomorphisms; complete
invariants: the supernatural number) \cite{jg60}.  This was eventually
superseded by the Elliott's classification of AF \cstar algebras
(building blocks: finite dimensional \cstar algebras; embeddings:
\star homomorphisms; complete invariants: dimension groups)
\cite{gae76}. 

Researchers in \cstar algebras have proceeded on to the study of
direct systems of \cstar algebras which are not finite dimensional;
researchers in non-self-adjoint algebras have discovered that for
systems of finite dimensional (non-self-adjoint) operator algebras the
limit algebras are highly dependent on both the choice of building
block algebra and the choice of admissible embeddings and that the
variety of limits which appear is exceptionally large.  
Given the
central role that triangular operator algebras play in the theory of
non-self-adjoint operator algebras, it is not surprising that the
best understood limits in that domain
are those for which the building block algebras
are drawn from the family of full upper triangular matrix algebras
($T_n$'s). 

Indeed, the first non-self-adjoint classification 
results were obtained
for direct limits of $T_n$'s with refinement and with standard
embeddings \cite{rlb90,ppw90,scp90b}.  This was soon followed by the
classification of alternation algebras \cite{hp92,ytp92}.

These three families are subsumed by the family of direct systems in
which the building blocks are $T_n$'s and the embeddings are order
preserving star extendible homomorphisms.  Order preservation, which
has its roots in \cite{ms89,scp92}, was analysed in great detail in
\cite{dh95} for the $T_n$ context (and for direct sums thereof); in
particular, a complete, albeit complicated, classification was given
for order preserving limits of full upper triangular matrix algebras. 

At a greater level of generality, triangular subalgebras of AF 
\cstar algebras, with the diagonal a canonical (Cartan) masa in the 
\cstar algebra, have played a major role in the theory of limit
algebras.  In particular, the spectrum (originally called the
fundamental relation) serves as a complete invariant for these TAF
algebras \cite{scp90a,scp92bk}.  A much broader class of 
similarly coordinatised  algebras is
obtained by permitting the building block algebras to be digraph
algebras (including nontriangular algebras).  Amongst these, a number
of subclasses have proven to be classifiable in terms of scaled
abelian groups.  A notable example is the family of $2n$-cycle
algebras, which are of particular interest since these are the
simplest operator algebras with non-trivial homology.  With rigid star
extendible embeddings, these have been classified in
\cite{dp99,dp97,scp99}. 

Like triangular algebras, nest algebras have played a critical role,
since its inception, in the theory of non-self-adjoint algebras.  This
paper is devoted to a study of direct systems in which the building
blocks are finite dimensional nest algebras.  For invariants, we turn
to the Grothendieck dimension distribution groups, which were first
developed in \cite{scp}.  These will be denoted by $\GG$.  This
invariant reduces to $K_0$ in the \cstar algebra context and to a dual
form of the spectrum in the TAF algebra context.  Dimension
distribution groups, augmented as necessary by various ancillary
structures, form an effective invariant for additional families of
systems -- as we shall see with finite dimensional nest algebras --
and help unify the classification of direct systems and of limit
algebras.  

It is, however, necessary to restrict the class of admissible
embeddings between finite dimensional nest algebras in order to
obtain classification theorems in these terms.  We shall show that
concepts of order preservation  and conservation give rise to diverse
embeddings which are tractable in this way.
  Order preservation for embeddings is defined in terms of
order preservation for isometries; in \cite{apd98}, Donsig extended
these concepts from the $T_n$ context to the digraph algebra context.
The definitions are made with respect to a self-adjoint subalgebra of
the original algebra; in the triangular case the diagonal is the only
reasonable choice for this subalgebra.  When $A$ is a finite
dimensional nest algebra, $A \cap A^*$ is a natural choice for this
subalgebra; we adopt Donsig's definition in the context of this paper.

We shall find, however, that critical preliminary theorems
($K_0$-uniqueness and $\GG$-lifting) are valid for broader classes of
embeddings and so classification can be achieved beyond the realm of
order preserving systems.  A related concept, order conservation, is
central for this extension.  The reason this new concept appears for
finite dimensional nest algebras, and 
did not appear for $T_n$'s, is that the
diagonal order on projections is not anti-symmetric in the former
context as it is in the latter.  With the understanding that we
consider only regular partial isometries (see section~\ref{S:regem}),
we define an order preserving partial isometry  to be
one for which conjugation by the partial isometry preserves the
diagonal order ($p \preceq q$ if, and only if, $vpv^* \preceq vqv^*$,
where $p$ and $q$ are subprojections of the initial space of the
normalizing partial isometry $v$).  We define \orcon partial
isometries to be ones which respect, but do not necessarily preserve,
the diagonal order (in the sense that $p \prec q$ and 
$vqv^* \prec vpv^*$ cannot simultaneously hold).  In the $T_n$
context, these two concepts coincide.  But for finite dimensional nest
algebras in which at least some of the atoms have rank greater than
one, the order preserving partial isometries form a proper subset of
the \orcon partial isometries.

Order preserving embeddings are defined to be those embeddings which
map order preserving partial isometries to order preserving partial
isometries.  Similarly, \orcon embeddings map order
conserving partial isometries to \orcon partial isometries.
In the $T_n$ context, these two families of embeddings coincide.  In
general, however, neither family contains the other.  Examples
illustrating these facts are given at the end of Section~\ref{S:opemb}. 

There are local variants of these concepts which are of importance.
An embedding is locally order preserving if it maps rank one (regular)
partial isometries to order preserving partial isometries 
and is locally order
conserving if it maps rank one (regular) partial isometries to order
conserving partial isometries.  Clearly, each property implies its
local variant; furthermore, local order preservation implies local
order conservation (since order preserving partial isometries are
\orcon).  The following  summarizes these relationships: 
\begin{align*}
OP \Longrightarrow LOP \Longrightarrow &LOC \\
OC \Longrightarrow &LOC
\end{align*}
Thus, the locally order conserving embeddings form the broadest
family, containing all the other families.  Happily, two properties
essential for dimension distribution group
 classification, $K_0$-uniqueness and $\GG$-lifting, are
valid for all locally order conserving embeddings.  
See Theorem~\ref{T:r1uni} and Theorem~\ref{T:lift}. 
 These relationships and facts
indicate the importance of a purview broader than order preservation. 

Based on these results, we prove that the 
dimension distribution group,
together with ancillary structures, is a complete invariant for
systems of finite dimensional nest algebras of four types, one for
each preservation/conservation property.  The classification is up to
regular isomorphism of the appropriate preservation/conservation
type. See Theorems~\ref{T:clas1}, \ref{T:clas2}, \ref{T:clas3},
 and \ref{T:clas4}.

The classification of operator algebra direct limits involves
substantial technical difficulties and has been deferred.
However, for algebraic direct limits (locally finite algebras),
we obtain two classification theorems -- one for limits of \orpre
systems using the dimension distribution group and an \orpre scale
and one for limits of \orcon systems using the dimension distribution
group and an \orcon scale.
See Theorems~\ref{T:occlas} and \ref{T:opclas}.  

We now briefly indicate related classification schemes for
nontriangular limit algebras.
For nest algebras with just two atoms, star extendible embeddings are
automatically \orcon (although not necessarily order
preserving). The corresponding $T_2$-limit algebras $A$ were classified
in  \cite{scp92} and in  \cite{dh97} in terms of an
augmentation of $K_0$ data by a binary relation on $\Sigma_0(A) \times
\Sigma_0(A)$ where $\Sigma_0(A)$ is the $K_0$ scale of $A$. Also,
Donsig \cite{apd98} has obtained a complete classification of order
preserving limits of general finite dimensional nest algebras (and more
general chordal algebras) by means of a ternary relation augmentation
of $K_0$ data.
On the other hand, away from order conservation,
  in  \cite{scp_b}  $T_r$-limit
algebras, with general star extendible embeddings, 
are classified in terms of metrized
dimension module invariants.

Here is a guide to the four remaining sections of this paper.
Section~\ref{S:regem} discusses regular embeddings between finite
dimensional nest algebras.  In particular, we show that local
regularity implies regularity in the finite dimensional nest algebra
context. This result is not valid for general digraph algebras.  

Section~\ref{S:opemb} introduces the various concepts of order
preservation and order conservation.  The $K_0$-uniqueness theorem for
locally order conserving embeddings is proven.  This is a stronger
result than we actually need for classification; $\GG$-uniqueness
would suffice.  This section also contains a digression (which the
reader interested solely in classification may skip) into the
structure of \orcon embeddings.  We obtain a result
(Corollary~\ref{C:struct}) analogous to the theorem in~\cite{dh95}
which states that an order preserving embedding in the $T_n$ context
is an ordered sum of refinement embeddings.  For this purpose, we need
to introduce refinement type embeddings (a generalization of
refinement embeddings) and $T_2$-degenerate embeddings, which are
trivially \orcon but may not be of refinement type.
Curiously, in order for the theorem to be valid, we must assume that
every atom of the domain algebra has rank two or greater.  Thus,
a structure therem for \orcon embeddings is available when
all atoms have rank one or when all atoms have rank greater than one.
For the case with mixed rank one and multi-rank atoms, we show by
counterexample  that the theorem is false.  (A study of examples
suggests that some sort of general structure theorem is valid, but
that the basic order irreducible types are sufficiently complicated
and special that such a theorem may be of dubious interest.)

Section~\ref{S:iso} studies isomorphism of the various systems of
interest.
To obtain classifications of the algebraic limits up to star extendible
isomorphism we require the result that a general star extendible
isomorphism between the limit algebras of \lorcon systems
 necessarily derives from a regular
isomorphism between the given  systems (Theorem~\ref{T:istmp}). 
This guarantees that our invriants for limit algebras
are well-defined.
Moreover,  we refine this result 
(in Theorem~\ref{T:autooc}) to enable the formulation of 
 invariants defined in terms of \orpre or
\orcon maps (rather than regular maps).

Section~\ref{S:invar} gives a brief introduction to partial isometry 
dimension distribution groups.  For more detail, the reader is
referred to~\cite{scp}.  The $\GG$-lifting theorem for locally order
conserving embeddings is proven (and yields $\GG$-lifting for all the
other types of embeddings) and various scales are defined for
dimension distribution groups.  The six classification theorems
described earlier in the introduction are proven and the section
concludes with several examples.  Readers may wish to glance at these
examples before perusing the whole paper.

The authors  thank Allan Donsig and Paul Haworth for
 several helpful comments.

\section{Regular and locally regular embeddings} \label{S:regem}

A finite dimensional nest algebra is a unital  operator
algebra $A$ in $M_n = M_n({\mathbb C})$
which has a block upper triangular form with respect to an ordered
decomposition of the identity as an orthogonal sum of
projections: $1 = q_1   + q_2   + ... + q_l$. 
Such an algebra $A$ consists of all those
operators $a$ in $M_n$  for which $(1 - f)a f = 0$ for all projections $f
= q_1  +...+ q_k$, $1 \le k \le l.$
The projections $f$, together
with the zero projection, form the invariant projection nest for $A$,   
while the projections   $q_i$ are the atomic interval
projections of A.

The nest algebra $A$ induces, through its atomic interval projections,
a block structure on $M_n$ whereby each matrix $v$ in $M_n$
can be viewed as an $l \times l$
block matrix $ v = (v_{ij})$ where $v_{ij} = q_ivq_j$.
The elements of the
nest algebra are precisely the block upper triangular matrices with
respect to this block structure.  Unless some other block structure is
explicitly specified, it will always be assumed that elements of a
nest algebra have the block structure natural to the nest.

\begin{definition} \label{D:reg}
A partial isometry $v$ in $M_n$  is said to be \emph{regular} with respect to
a block structure if each  block matrix entry  is a partial isometry.
\end{definition}

Note that any two block matrix entries  in the same row 
of a regular partial isometry have orthogonal final
spaces and any two block matrices in the same column have orthogonal
initial spaces. 

\begin{definition} \label{D:std}
A partial isometry $v$ in a finite dimensional C*-algebra
is said to be 
 \emph{standard}, with respect
to a matrix unit system if it is a sum of some of these matrix units.
\end{definition}

Note that if a partial isometry is standard with respect to a 
matrix unit system, then it is regular with respect to any
block structure compatible with that matrix unit system.

An algebra homomorphism      $\phi   \colon A_1 \to     A_2$  
between  two   nonself-adjoint operator algebras
is said to be star extendible if it 
is the restriction of a  star algebra homomorphism
$\tilde{\phi} $ between the generated C*-algebras.
Plainly, $\tilde{\phi}$ is uniquely determined by $\phi$.
In the next definition and throughout the paper
we confine attention to such 
star-extendible homomorphisms.

\begin{definition} \label{D:rlr}
Let $\phi \colon A_1 \to A_2$ be 
a star-extendible homomorphism between 
finite
dimensional nest algebras.  Then  $\phi$ is said to be
\emph{regular} if $\phi$ is a direct sum of multiplicity one
embeddings. Also, $\phi$ is {\em locally regular\/} if
$\phi(v)$ is a regular partial isometry in $A_2$ whenever $v$ is
a regular partial isometry in $A_1$.
\end{definition}

\begin{remarks}
If $\phi$ is a regular embedding then $\phi$ is locally
regular.  We prove the converse in Theorem~\ref{T:reg}
below.  Although the definition above makes sense for 
maps between digraph algebras
it is not generally true that
local regularity is equivalent to regularity.
For example, consider the algebra $E \subseteq T_3$
which is spanned
 by all the matrix
units of $T_3$ except $e_{23}$. Then one can verify that every
star extendible map $\phi \colon E \to E \otimes M_n$ for $n \ge 2$ 
is locally regular, and yet there are irregular (i.e.\ nonregular)
embeddings.
See also the 4-cycle algebra example, Example 2.2, in \cite{dp99}.

In order to show that a star-extendible map
$\phi$ is locally regular, it is sufficient to
show that $\phi(v)$ is a regular partial isometry 
in $A_2$ whenever $v$ is a
rank-one regular partial isometry in $A_1$.  To see this,
let $\tilde{Q}_1, \dots , \tilde{Q}_m$ be the atomic interval
projections of $A_1$ and $Q_1, \dots, Q_n$ the atomic interval
projections of $A_2$.  Let $v$ be a regular partial isometry 
in $A_1$. Then $v_{ij} = \tilde{Q}_i v\tilde{Q}_j$
is a partial isometry, for all $i$ and $j$.
Each $v_{ij}$, in
turn, can be written a sum of rank one partial isometries each of
which has the form $w = \tilde{Q}_i w \tilde{Q}_j$.  Thus, we may
write $v$ as a sum, $v_1 + \dots + v_p$ say, of rank one regular
partial isometries

Since $\phi$ is star extendible, $\phi(v) = \phi(v_1) + \dots
 + \phi(v_p)$ is a partial isometry and, by hypothesis, each
 $\phi(v_i)$ is a regular partial isometry in $A_2$.  From this
we may conclude that the set of initial projections of the
$\phi(v_i)$ are pairwise orthogonal, the set of final projections
of the $\phi(v_i)$ are pairwise orthogonal, and both the initial 
projection and the final projection of each $\phi(v_i)$ commute
with each atom $Q_t$.  The last statement follows from the fact
that every $Q_s \phi(v_i) Q_t$ is a partial isometry.

It now follows that, for all $s$ and $t$, the set of initial 
projections of the elements 
$Q_s \phi(v_i) Q_t$ $(i = 1, \dots p)$ are pairwise orthogonal;
the same is true for the set of final projections.  Consequently,
$Q_s \phi(v) Q_t = Q_s \phi(v_1) Q_t + \dots + Q_s \phi(v_p) Q_t$
is a partial isometry.  Since this is true for all $s$ and $t$, we see
that $\phi(v)$ is regular in $A_2$ and that $\phi$ is locally regular.

We may, in fact, go a bit further than reducing local regularity 
to the action of $\phi$ on rank-one regular partial isometries.
If $\{e_{ij}\}$ is a matrix unit system in $A_1$ compatible
with the block matrix structure of $A_1$ and if $\phi(e_{ij})$
is regular in $A_2$ for all $e_{ij}$, then $\phi$ is locally
regular.

If $\phi$ is locally regular
and $\{e_{ij}\}$ is a matrix unit system compatible with $A_1$, then
for each matrix unit $e_{ij}$ there is a matrix unit system for 
$A_2$ with respect to which $\phi(e_{ij})$ is a sum of matrix units;
i.e.\ $\phi(e_{ij})$ is standard.  If, moreover,  
there is a single matrix unit
system for $A_2$ with respect to which all $\phi(e_{ij})$ are
standard, then it follows that $\phi$ is regular.  
Obtaining such a system will be the strategy for the proof
of the following theorem.

\end{remarks}

\begin{thm} \label{T:reg}
Let $\phi   \colon A_1 \to     A_2$
be a star extendible homomorphism between finite 
dimensional nest algebras.
If $\phi$ is locally regular then $\phi$ is regular.
\end{thm}

\begin{proof}
The nest algebra $A_2$ and the \cstar algebra it generates have
an $r \times r$ block structure induced by the $r$ atoms in the nest
for $A_2$.  In this block structure, the diagonal blocks are square
matrices of varying size and the off-diagonal blocks are rectangular.
Elements of $A_2$ and $C^*(A_2)$ will be 
written as $r \times r$ matrices, $v = (v_{st})$, $1 \leq s,t \leq
r$, where each $v_{st}$ is a matrix of appropriate size.

We may assume that $A_1 = T_n$.
Consider the star algebra extension of $\phi$. This is an algebra
injection $\phi \colon M_n \to C^*(A_2)$ 
such that the image of each standard matrix unit of $M_n$
is a regular partial isometry. 
We wish to show that $\phi $ is regular with respect to the 
$r$ by $r$ block structure. 
By the remarks above this is equivalent to
showing the following with $v^{ij} = \phi(e_{ij})$.

Let $v^{ij} = (v^{ij}_{st}), $with $ 
1 \le s,t \le r, 1 \le i,j \le n,$ be a matrix unit system
for a subalgebra of $C^*(A_2)$ which is isomorphic to $M_n$.
Suppose moreover that each $v^{ij}$ is regular with respect to the
$r$ by $r$ block structure, so that each block matrix entry 
$v^{ij}_{st}$    of
each $v^{ij}$ is itself a partial isometry.
We show that there is a matrix unit 
system for $C^*(A_2)$, which 
consists of (rank one) regular partial isometries, such that each
$v^{ij}$ is a sum of some of these matrix units.

To see this consider first a
product $v^{ij}v^{jk}$. The ~$1,1$~ block entry is
given by the sum
\[
v_{11}^{ij}v_{11}^{jk} + 
v_{12}^{ij}v_{21}^{jk} +
 \dots +
v_{1r}^{ij}v_{r1}^{jk}.
\]
Since  $v^{ij}$ is regular, 
 the
partial isometries \,$v_{11}^{ij}, \dots , v_{1r}^{ij}$ \,have orthogonal range
projections and so 
the operators of the sum 
have orthogonal range projections.
For similar reasons the domain projections are pairwise orthogonal.
(We are not assuming here that these products are partial isometries.)
Since, by hypothesis, the product $v^{ij}v^{jk}$ is a regular
partial isometry, it follows that the sum above is a partial isometry,
and therefore, by the orthogonality of domain and range
projections, each of
the individual products
\[
v_{11}^{ij}v_{11}^{jk}, \,v_{12}^{ij}v_{21}^{jk}, \dots , 
\,v_{1r}^{ij}v_{r1}^{jk}
\]
is a partial isometry.

Now, since, for example, $v_{11}^{ij}v_{11}^{jk}$ is a partial isometry
it follows that the range projection of $v_{11}^{jk}$ commutes with the
domain projection of $v_{11}^{ij}$. Abusing notation somewhat,
and regarding the entry operators $v^{ij}_{st}$ as identified with
operators in $C^*(A_2)$, it follows, by considering other block entries,
that for all 
$i,j,k,l,s,t,u,v$ the range projection of $v^{ij}_{st}$ commutes with the
domain projection of $v^{kl}_{uv}$. Note also that the domain
projections and the range projections commute amongst themselves.
Furthermore it is clear
that these projections commute with the 
projections in the centre of the block diagonal subalgebra of
$C^*(A_2)$.
Choose now a maximal family of rank one projections $p$ which commute
with all these projections and are dominated by $v_{11}^*v_{11}$.
Then for each such projection  the set of operators $w_{ij} = 
v^{i1}pv^{1j}$ satisfy the relations of a matrix unit system. The
projections $P = \Sigma_i  v^{i1}pv^{1i}$ are pairwise orthogonal and
decompose $\phi$ as a direct sum of multiplicity one embeddings, as
desired.
\end{proof}

\section{\Orcon  and \orpre  embeddings} \label{S:opemb}

In this section  we define 
\orcon and \lorcon  star extendible
embeddings of finite dimensional nest algebras.  We show 
that these embeddings are classified
 up to inner conjugacy by $K_0$ and we make use of this in the
 classifications of Section~\ref{S:invar}.  This 
inner conjugacy classification can also 
be viewed as a classification 
of those subalgebras $A_1$  of a fixed nest algebra $A_2$
where $A_1$ is a nest algebra in $C^*(A_1)$  and where
the inclusion of $A_1$  is an \orcon inclusion.  

In \cite{dh95}, it is shown that unital order preserving embeddings in 
the $T_n$ context are direct sums of refinement embeddings.  We will
extend this result (slightly) to  non-unital embeddings and will
obtain an analogous decomposition for \orcon embeddings
between nest algebras, provided that all atoms in the domain algebra
have rank greater than one.  The decomposition of $\phi$ will be into
an ordered sum of order irreducible embeddings of two basic types to
be defined below.

In this paper we consider a general finite dimensional operator
algebra $A$ to be a nest algebra if it is star extendibly isomorphic
to a nest algebra. This means that the generated \cstar-algebra 
$C^*(A)$
is isomorphic to $M_p$ for some $p$ and $A$ is a nest algebra in
$C^*(A)$.  \medskip

\begin{defn}  \label{D:pi}
Let  $A \subseteq    M_n$   be a finite dimensional nest algebra
with atomic interval projections  $q_1, q_2, \dots , q_l$. 
\begin{enumerate}
\item 
 The set
$\Pisom_{\reg}( A)$  is the set of regular partial isometries
with respect to the block structure from $q_1, \dots, q_l$.

\item
 The set $\Pisom_{\oc}( A)$   denotes the subset 
of $\Pisom_{\reg}( A)$ consisting of {\em \orcon\/} 
partial isometries $v$; that is, those for which
\[
q_ivq_j \ne 0 \implies q_svq_t = 0, \quad \text{for } s < i \quad  
\text{and} \quad t > j.
\]

\item
A star extendible embedding between finite dimensional
nest algebras is {\em \orcon\/} if it 
maps \orcon partial
isometries to \orcon partial isometries.

\item
 An embedding $\phi \colon A_1 \to A_2$ is 
{\em \lorcon\/} if $\phi(v) \in \Pisom_{\oc}(A_2)$
for each rank one element $v$ in $\Pisom_{\reg}(A_1)$. 
\end{enumerate}
\end{defn}

If $v \in \Pisom_{\reg}(A) $ and $\rank v = 1 $ then
$v \in \Pisom_{\oc}(A) $.  It follows from 
this observation, the remarks 
following Definition~\ref{D:rlr}, and Theorem~\ref{T:reg} that a
\lorcon embedding between finite dimensional nest
algebras is automatically a regular embedding.

In particular, \orcon star extendible embeddings
are regular star extendible embeddings and so admit 
decompositions as direct sums of multiplicity one embeddings.

For an explicit example of an \orcon  embedding
consider the map $\phi$ from $T(2,2,2)$ to $T(6,8,10)$
given by 
\[
\begin{bmatrix}
a&x&z\\
 &b&y\\
 & &c
\end{bmatrix}
\longrightarrow 
\left[
\begin{array} {c c c | c c c c |  c c c c c}
a& &x& &z& & & & & & &\\
 &a& &x& &z& & & & & &\\
 & &b& &y& & & & & & &\\
\hhline {------------}
 & & &b& &y& & & & & &\\
 & & & &c& & & & & & &\\
 & & & & &c& & & & & &\\
 & & & & & &a& &x& &z&\\
\hhline {------------}
 & & & & & & &a& &x& &z\\
 & & & & & & & &b& &y&\\
 & & & & & & & & &b& &y\\
 & & & & & & & & & &c& \\
 & & & & & & & & & & &c \\ 
\end{array}
\right].
\]
This map has an order irreducible {\it ordered} sum
decomposition $\phi_1 + \phi_2 + \phi_3$
where $\phi_1$ is a multiplicity 2 map, with range dominated by
the sum of the first two atoms of the codomain algebra,
and $\phi_2$ and $\phi_3$ are multiplicity one embeddings.
The summands here are all of $T_2$-degenerate type 
(in the sense below).

In general, the \orcon  partial isometries in $A$ are those 
regular partial isometries whose
block matrix supports have a staircase form, such as
\[
\begin{bmatrix}
 * & * & * & 0 & 0 & 0\\
  & 0 & * & 0 & 0 & 0\\ 
  &   & * & * & * & 0\\ 
  &   &   & 0 & * & 0\\
  &   &   &   & * & *\\
  &   &   &   &   & *
\end{bmatrix}. 
\]

\begin{thm} \label{T:r1uni}
Let $A_1 $ be a finite dimensional
nest algebra.  Let $\phi $ and $\psi $ be two
\lorcon regular embeddings of $A_1 $ into a second
finite dimensional nest algebra, $A_2$.  Then $\phi $ and $\psi $ are
inner conjugate  if, and only if, $K_0(\phi) = K_0(\psi) $.
\end{thm}

\begin{proof}
We first prove the theorem for the special case in which every atom
of $A_1$ has rank one (i.e.\ $A_1$ is some $T_n$).
Select a matrix unit system for $C^*(A_1)$ compatible with $A_1$; for
each pair of (rank one) atoms $s$ and $r$, let $m(s,r)$ denote the
matrix unit with initial projection $r$ and final projection $s$.
Note that
$K_0(\phi) $ is a matrix whose columns are indexed by the atoms of
$A_1 $ and whose rows are indexed by the atoms of $A_2 $.  The entry
in the column indexed by the atom $r$ of $A_1$ and
the row indexed by the
atom $Q$ of $A_2$ is the non-negative integer ~$\rank \phi(r)Q $.
We assume the order of the rows and columns of this matrix
reflect the usual $\prec$ ordering of the atoms of $A_1$ and $A_2$.

Let $r$ be an atom of $A_1$ with the property that the first non-zero
entry in the column of $K_0(\phi)$ 
indexed by $r$ is less than or equal to the first
non-zero entry in any other column.  Let $Q$ be the atom which indexes 
the row containing the first non zero entry in column $r$.  Now assume 
that $s$ is some other atom of $A_1$.  First, assume $s \prec r $.
Let $P$ be the atom of $A_2$ which gives the first non-zero entry in
column $s$.  Let $E_r = \phi(r)Q$ and $F_s = \phi(s)P $.  The fact
that $\phi$ is \lorcon  ensures  that
$\phi(m(s,r)) E_r \phi(m(s,r))^* \leq F_s $. 
To see this note that for 
an appropriate choice of matrix units in $A_2$,
 $\phi(m(s,r))$ induces an \orcon
  bijection from the rank one matrix unit
subprojections of $\phi(r) $ to the rank one matrix unit
subprojections of $\phi(s)$.  Since $E_r$ and $F_s$ correspond to the
first block occurences of these matrix unit projections and since, by
the choice of $r$, $F_s$ has at least as many of them as $E_r$ does,
it follows that $\phi(m(s,r))$ conjugates $E_r$ to a subprojection of
$F_s$.  

Let $E_s = \phi(m(s,r)) E_r \phi(m(s,r))^*$.  If  $r \prec s $, we argue in 
a similar fashion to obtain a projection $E_s$ which is a
subprojection of the index atom for the row containing the first
non-zero entry in column $s$ and which satisfies
$E_r = \phi(m(r,s)) E_s \phi(m(r,s))^* $.  The fact that $\phi$ is an
embedding implies that, for any two atoms $s \prec s'$, 
$E_s = \phi(m(s,s') E_{s'} \phi(m(s,s'))^*$.

Note that if $s$ and $s'$ are distinct atoms of $A_1$, then the
projections $E_s $ and $E_{s'} $ are orthogonal.
This is true even if $E_s $ and 
$E_{s'} $ are subprojections of the same atom of $A_2 $, since 
$E_s \leq \phi(s)$, $E_{s'} \leq \phi(s')$ and $\phi(s)$ and
$\phi(s') $ are orthogonal projections.

Let $ P_1 = \sum_s E_s $, where the sum runs over all atoms of
$A_1$. Let $\phi_1 = P_1 \phi P_1$.  It follows from the way that 
$\phi_1 $ is defined that $\phi_1 $ and $(1-P_1) \phi (1-P_1) $
are \lorcon 
 embeddings of $A_1$ into $A_2$ and that
$\phi = \phi_1 + (1-P_1) \phi (1-P_1)$. The reason that $\phi_1$ and  
$(1-P_1) \phi (1-P_1) $ are \lorcon is that any
subpartial isometry of an \orcon  partial isometry which is
obtained by left and right multiplication by a block diagonal
projection is again \orcon.    Note that
$K_0((1-P_1) \phi (1-P_1)) $ is obtained from $K_0(\phi) $ by
subtracting from the first non-zero entry in each column of
$K_0(\phi) $ the integer in column $r$, row $Q$.

Since $K_0(\psi) = K_0(\phi) $, we may construct a projection 
$P^{\psi}_1 $ which has the same properties with respect to $\psi$
that $P_1$ has with respect to $\phi$.  $P^{\psi}_1 $ is a sum of
projections, one for each atom of $A_1$; the projection in the sum
associated with an atom $s$ has the same rank as $E_s$ and is a
subprojection of the same atom of $A_2$ as $E_s$.  Consequently, there 
is a  unitary element of the diagonal $A_2 \cap A_2^*$ such that if
$\tilde{\psi} = u^* \psi u $, then 
$P_1 \phi P_1 = P_1 \tilde{\psi} P_1 $ and
$(1-P_1) \tilde{\psi} (1-P_1) $ is an embedding.

We may now repeat this procedure with the two embeddings,
$(1-P_1) \phi (1-P_1) $ and $(1-P_1) \tilde{\psi} (1-P_1) $. When we
do so, we obtain a projection $P_2$ which is orthogonal to $P_1$ and
an embedding $\hat{\psi}$ which is inner conjugate to $\psi$ such that
$P_2 \phi P_2 = P_2 \hat{\psi} P_2 $. Furthermore, because we are
really working in essence on $1-P_1$, we can select the unitary which
implements the conjugacy so that it is the identity on $P_1$.  With
this done, we also have $P_1 \phi P_1 = P_1 \hat{\psi} P_1 $ and hence 
$(P_1+P_2) \phi (P_1+P_2) = (P_1+P_2) \hat{\psi} (P_1+P_2) $.  It is
now clear that if we continue in this fashion, after finitely many
steps we obtain an inner conjugacy between $\phi$ and $\psi$.

Slight modifications of the argument given above yield the proposition 
when $A_1$ is a general finite dimensional nest algebra.
Alternatively, we can deduce the general version of the proposition
from the special case above by the following 
 `general principles' technique.
If $A_1$ is a general nest algebra, select a subalgebra $B$ of $A_1$
which is isomorphic to some $T_n$  such that each block
of $A_1 \cap A_1^*$ corresponding to an atom of $A_1$ contains exactly
one rank one diagonal projection from $B$.
Note that $A_1$ is the algebra
generated by $B$ and $A_1 \cap A_1^* $. 
Let $\phi$ and $\psi$ be two \lorcon\ embeddings such
that $K_0(\phi) = K_0(\psi) $.  Since $\phi$ and $\psi$ are star
extendible, they are determined (up to inner equivalence)
 by their restrictions to $B$.  But 
$K_0(\phi |_B) = K_0(\phi) = K_0(\psi) = K_0(\psi |_B)$.  By the
argument above, $\phi |_B$ and $\psi |_B$ are inner equivalent; it
follows immediately that $\phi$ and $\psi$ are inner equivalent.
\end{proof}

\subsection{Structure theorems} \label{SS:str}

We now determine the way in which \orcon
embeddings decompose as ordered  
sums of embeddings of special types. Our original
 motivation for studying such decompositions was to use these
as a vehicle for obtaining $K_0$-uniqueness.  Theorem~\ref{T:r1uni},
however, provides $K_0$-uniqueness more directly and 
for a broader class of embeddings.  While the original motivation is
no longer germane, the structure theorem is  of independent
interest.  We will also comment on
connections with $K_0$-uniqueness, since these connections are
illuminating. 

In the next definition we write   $\prec_A$
to indicate the total ordering on the atomic projections of a nest algebra
$A$.

\begin{defn} \label{D:rftyp}
 (i)  A (possibly nonunital) star extendible embedding
$\phi \colon T_n \to T_m$ 
is said to be a {\it refinement type} embedding if $\phi$ is 
\lorcon  and if,
whenever  $e_1 \prec_{T_n} e_2$     and $e_1 \ne e_2$,
it follows that $f_1  \prec_{T_m} f_2$ for $f_1 \le \phi(e_1)$ and
$f_2 \le \phi(e_2)$.

(ii)  A (possibly nonunital) star extendible embedding 
$\phi \colon 
A_1 \to A_2$ between finite dimensional nest algebras is said to be of 
{\it refinement type} 
if there are maximal triangular subalgebras  $T_n \subseteq A_1$,
$T_m \subseteq A_2$
such that   $\phi(T_n) \subseteq T_m$ and the restriction map
$\phi \colon T_n \to T_m$
is of refinement type.
\end{defn}

A refinement type embedding is necessarily \orcon. To see this 
let $\rho \colon T_n \to T_m$ be a refinement 
embedding with star extension
$\rho \colon M_n \to M_m$ and let $A_1$ be a nest algebra containing
$T_n$. The algebra $A_1$ necessarily has a standard block upper
triangular structure. There is a smallest nest algebra $A_2$ containing 
$T_m$ and $\rho(A_1)$ and this algebra
 also has standard form.
Note that the images of the atoms of rank greater than one
are atoms of $A_2$ whilst the other atoms of $A_2$ are of rank one.
The map $\rho \colon A_1 \to A_2$ is clearly \orcon.
If $A_3$ is a nest algebra in $M_m$ containing $A_2$ then since
the block structure of $A_2$ refines that of $A_3$  it follows
that $\Pisom_{\oc} (A_2) \subseteq \Pisom_{\oc} (A_3)$. Thus
$\rho \colon A_1 \to A_3$ is \orcon and
such a map $\rho$ is a typical unital refinement type embedding.

A general refinement type embedding $\phi \colon A_1 \to A_4$ 
has the form $\psi \circ \rho$ where $\psi \colon A_3 \to A_4$
is a multiplicity one inclusion, and so is also \orcon.

Note that the argument above implies that the image of an atom
of rank greater than one, under a refinement type embedding,
is dominated by an atom of the range.

\begin{defn} \label{D:orsum}
Let  
$\phi_i \colon A_1 \to A_2$, $1 \le i \le r$,
be star extendible algebra homomorphisms between finite dimensional nest
algebras $A_1   , A_2$    such that the projections     $\phi_i(1)$
are pairwise
orthogonal.  Then the star extendible embedding  
$\phi_1 + \dots + \phi_r$ is said to be an {\it ordered sum }
if \[
q\phi_i(1) \ne 0 \implies q'\phi_j(1) = 0 
\quad \text{for} \quad q \prec q', j < i
\]
for all atomic interval projections $q,q'$ of $A_2$, and for all $i$.
\end{defn}

What this means can also be described in terms of a choice of a maximal set 
${\cal P}
=  \{e_1  ,e_2  ,..., e_m\}$        of rank one projections in $A_2$
whose indexing
is consistent with the partial ordering  $\prec_{A_2}$
       on the set  ${\cal Q}$ of
atomic interval projections $q$ (in the sense that if 
$e_i \le q$, $ e_j \le q', i \le j$, 
then $q \preceq_{A_2} q'$) and is such that each projection
$\phi_i(1)$ is a sum of some of the projections of ${\cal P}$.
With any such choice of ${\cal P}$
the condition of the last definition holds if the ${\cal P}$-support
of the projections $\phi_1(1),\phi_2(1), \dots ,\phi_r(1)$
appear in order, that is,
\[
\max \{j : e_j\phi_i(1) \ne 0\} < \min \{j : e_j\phi_{i+1}(1) \ne 0\}
\]
for each  $i = 1,...,r.$

\begin{defn} \label{D:orred}
Let $\phi \colon A_1 \to A_2$    be a star extendible algebra
homomorphism with ordered decomposition
$\phi_1 + \dots + \phi_r$.
Such a decomposition is said to be \emph{order reducible}   if for some
index  $i$
the map $\phi_i$
has a non-trivial ordered sum decomposition.
\end{defn}

The following proposition is evident.

\begin{prop} \label{P:ordec}
 A star extendible embedding between finite dimensional 
nest algebras
admits an order irreducible ordered sum decomposition which is unique up
to a permutation of multiplicity one summands which map into the same
diagonal block. Moreover the embedding is \orcon  if, and only
if, each summand is an \orcon  embedding.
\end{prop}

This  proposition provides a useful perspective and to
illustrate this we prove a key result from  \cite{dh95}
generalised to the case of nonunital maps.

\begin{prop} \label{P:embfo}
Let $n \ge 3$ and let    $\phi  \colon T_n  \to    T_m$    be a 
star extendible embedding.
Then  $\phi$ is an \orcon embedding if, and only if, $\phi$ is an
ordered sum  $\phi = \phi_1 + \dots + \phi_p$, 
where each  $\phi_j$ is a refinement type embedding.
\end{prop}

\begin{pf}
Let  $\phi$ be \orcon with order-irreducible ordered
decomposition 
\[
\phi_1 + \phi_2 + \dots + \phi_l
\]
 and let  ${\cal P} = \{e_1, \dots ,e_n\}$
be the atomic interval projections of $T_n$, listed in their
natural order.  Let 
$Q_1, Q_2, \dots ,Q_j$
be consecutive interval projections of $T_m$ (consecutive in the 
sense of there being no gaps)  chosen with
the following properties: 
\begin{align*}
  &Q_i\phi(e_i) = Q_i\phi(1), \text{ for each $i$}, \\
  &\rank (Q_1\phi(e_1)) = \rank (Q_2\phi(e_2)) = \dots 
    = \rank (Q_{j-1}\phi(e_{j-1})) \ge \rank (Q_j\phi(e_j)), \\
 & \text{$j$ is maximal},
\end{align*}
 where either $j = n$ and the last inequality is 
an equality, or $j < n $ and
the last inequality is strict.    In the first case,
since  $\phi$ is \orcon,
 it follows that $\phi$  has an ordered
decomposition    $\phi_1' + \phi_2'$, where $\phi_1'$ is of
refinement type; 
 thus, by order irreducibility, 
$\phi_1 = \phi_1'$. We show that the second case 
cannot hold.  This, and induction,
completes the proof of the proposition.

Let us illustrate this situation with the following diagram,
where $j=4$, $a_1, a_2, a_3$ are the rank one diagonal
projections of $T_m$, in order, comprising $Q_1 \phi(e_1)$;
$b_1, b_2, b_3$, those of $Q_2 \phi(e_2)$; $c_1, c_2, c_3$,
those of $Q_3 \phi(e_3)$; and $d_1, d_2$, those of
$Q_4 \phi(e_4)$:
\[
\underbrace{a_1, a_2, a_3}_{Q_1}, \underbrace{b_1, b_2, b_3}_{Q_2},
\underbrace{c_1, c_2, c_3}_{Q_3}, \underbrace{d_1, d_2}_{Q_4},
x, \dots, y, \dots .
\]

Now let $x$ denote the first rank one projection of $T_m$ which
follows $Q_j$ and satisfies $x \phi(1) \neq 0$.  Since $\phi$ 
maps $T_n$ into $T_m$,
 the assumption on $j$ implies that $x$ is necessarily a
subprojection of $\phi(e_1)$ (and in our illustration could be
written as $a_4$).  Also, let $y$ be the first rank one projection
of $T_m$ such that $y \phi(e_j) \neq 0$ and $y Q_j = 0$.  
Such a projection exists
since $\phi$ is star extendible and so the projections
$\phi(e_1), \dots, \phi(e_j)$ have equal rank.  (In our illustration,
$y$ could be written as $d_3$.)

We can now see that order conservation is contradicted.  For note that
$\phi(e_{j-1,j})$ maps $y$ into $Q_{j-1}$.  Thus, while the partial
isometry $v = e_1 + e_{j-1,j}$ is \orcon, we have
$x \prec_{A_2} y$ and $\phi(v) y \phi(v)^* \prec_{A_2}
\phi(v) x \phi(v)^*$.
\end{pf}

Refinement type embeddings    $  \phi  \colon A_1 \to     A_2$
are determined, uniquely,  by their triangular 
restrictions,  $\phi_t \colon  A_1^t  \to  A_2^t$,  where
$A_1^t \subseteq  A_1$   is a maximal triangular subalgebra 
(isomorphic to $T_n$  for some  $n$)
    and where  $ A_2^t \subseteq  A_2$ 
is a similar appropriately
chosen triangular subalgebra. Compositions of refinement type
embeddings between triangular nest algebras are of refinement type; it
follows from this and the remark after Definition~\ref{D:rftyp}
that general refinement type embeddings are closed under
compositions.

Furthermore, refinement type embeddings between triangular 
algebras are determined up to inner conjugacy by
their  restrictions to the diagonal subalgebra.  Accordingly, it
follows (even without the use of Theorem~\ref{T:r1uni})
that refinement type embeddings  $\phi, \psi \colon A_1 \to  A_2$  
are inner conjugate if, and only if, $K_0\phi = K_0\psi$. 

\begin{defn}  \label{D:degty}
A star extendible embedding  $\phi \colon A_1 \to A_2$
between finite dimensional nest algebras
is said to be of $T_2$-\emph{degenerate type} if there 
exist atomic interval projections
$Q, Q'$ for $A_2$ such that $\phi(1) \le Q + Q'$.
\end{defn}

For an example of an embedding of $T_2$-degenerate type,
consider the multiplicity two embedding $\phi \colon T_4 
\otimes M_m \to T_2 \otimes M_n$ for which
\[
\left[ \begin{array}{cccc}
a&x&.&.\\
&b&y&.\\
&&c&z\\
&&&d
\end{array} \right]
\longrightarrow
\left[ \begin{array}{cccc|cccc}
a& & & &x&*&*&\\
 &a&x&.& & & &.\\
 & &b&y& & & &.\\
 & & &c& & & &z\\
\hhline{--------}
 & & & &b&y&*& \\
 & & & & &c&z& \\
 & & & & & &d& \\
 & & & & & & &d
\end{array} \right]
\]

Note that $\phi = \theta_1 + \theta_2$ where $\theta_1 , \theta_2$ are
rank one embeddings determined, respectively, by the block groupings
$(a, bcd)$ and $(abc, d)$.  The presence of $b$ in block two and $c$
in block one show that $\phi$ is not of refinement type.  Also $\phi$
is order-irreducible.  Since it can be checked that all order
irreducible embeddings from $T_3 \otimes M_n$ to $T_2 \otimes M_n$ are
of refinement type, $\phi$ can be viewed as the simplest order
irreducible non refinement type \orcon  embedding.

It can be shown, as in \cite{dh97}, that a $T_2$-degenerate
embedding is automatically regular.  Furthermore, every 
$T_2$-degenerate embedding is, essentially by default, \orcon ,
and so
Theorem~\ref{T:r1uni} implies that $T_2$-degenerate
embeddings have the $K_0$-uniqueness property.  Since it is
illuminating, we  give a direct proof of this.

Let $ q_1, \dots , q_p$ (resp. $Q_1, \dots , Q_n$)
be the atomic interval projections of $A_1$ (resp. $A_2)$.
Then, for a $T_2$-degenerate embedding 
$\phi$, one sees that $K_0\phi$ is an
$n \times p$ matrix with at most two nonzero rows, yielding a submatrix
\[
 \left[ \begin{array}{cccc}          
u_1 & u_2 & \dots & u_p\\
v_1 & v_2 & \dots & v_p
\end{array} \right]
\]
for which $v_i \le v_j$ and $u_i \ge u_j$, for $i \le j$,
and $u_i + v_i = r$, for all $i$, where $r$ is the multiplicity of $\phi$.

 The map $\phi \colon A_1 \to A_2$
admits a direct sum decomposition into rank one embeddings
belonging to  $p+1$ distinct inner equivalence classes.
The corresponding $K_0$ matrices for embeddings in these classes
have the form
\begin{align*}
  & \begin{bmatrix}          
1 & 1 & \dots & 1 & 1 & 1 \\
0 & 0 & \dots & 0 & 0 & 0
\end{bmatrix}, \\
  & \begin{bmatrix} 
1 & 1 & \dots & 1 & 1 & 0 \\
0 & 0 & \dots & 0 & 0 & 1
\end{bmatrix}, \\ 
  & \begin{bmatrix}
1 & 1 & \dots & 1 & 0 &0 \\
0 & 0 & \dots & 0 & 1 &1
\end{bmatrix}, \\
  & \qquad \qquad \vdots \\ 
 & \begin{bmatrix} 
0 & 0 & \dots & 0 & 0 & 0 \\
1 & 1 & \dots & 1 & 1 & 1
\end{bmatrix}. 
\end{align*}
If there are $r_k$ repetitions of the k$^{th}$ class 
then we say that $\phi$ has
multiplicity signature $\{r_1, \dots , r_{n+1}\}$. In this case
\[
K_0\phi =
 \left[ \begin{array}{cccc}
r_1 + r_2 + \dots r_p & r_1 + r_2 + \dots r_{p-1} & \dots & r_1\\
r_{p+1}  & r_p + r_{p+1}  & \dots & r_2 + r_3 + \dots r_{p+1}  
\end{array} \right],
\]
and so it follows that $K_0\phi$ determines the multiplicity
signature and thus the inner equivalence class of $\phi$.

When the following theorem is combined with Proposition~\ref{P:ordec}, 
it yields a structure theorem analogous to Proposition~\ref{P:embfo}
for \orcon embeddings defined on a nest algebra with no rank 
one atoms.

\begin{thm} \label{T:emfor}
 Let $\phi \colon  A_1 \to A_2$ be a regular star extendible
embedding between finite dimensional nest algebras which is \orcon
 and order irreducible, and suppose that $A_1$ has no
one dimensional atoms.  Then $\phi$ is either $T_2$-degenerate or is
of refinement type.  
\end{thm}

\begin{proof}
  Assume that $\phi$ is not $T_2$-degenerate.
Let $A_1$ have atomic interval projections $q_1, q_2, \ldots, q_n$ and
let $A_2$ have atomic interval projections $Q_1, Q_2, \ldots, Q_m$.
The hypothesis for $A_1$ implies that there exists a matrix unit
system $v^p_{ij}$, $1 \leq i \leq j \leq n$, $p = 0,1$,
for a copy of $T_n \otimes {\mathbb C}^2$ in $A_1$ such that $v^p_{ij}$ is
supported in  $q_i A_1 q_j,$ the $ (i,j)$ block subspace of $A_1$.
In particular note that for $i<n$ there are 
rank two \orcon  partial
isometries in $A_1$ of the form $v = v^0_{i,j_1} + v^1_{i, j_2}$.
Similarly, there are rank two \orcon  partial
isometries supported in a single block column.

Let $\phi = \theta_1 + \cdots + \theta_r$ be the multiplicity one
decomposition of $\phi$, with $r >1.$  
Note first that since $\phi$ is \orcon,
any subsum, such as $ \theta_s + \theta_t$, is also
\orcon.  Secondly, observe that for each $s$,
the partial isometry $  \theta_s
(v^0_{1n})$ is supported in a single off-diagonal block.  For if not,
then $Q_i \theta_s (v^0_{1n})Q_i = \theta_s (v^0_{1n})$ for some $s$ and
$i$, and so the range of $\theta_s$ lies in $Q_i A_2 Q_i$.  On the other
hand $\phi$ is \orcon  and it follows readily that $\theta_s$
must be an order summand of $\phi$, and hence equal to $\phi$ by 
order irreducibility, contrary to the assumption $r >1$. 

Consider an index $t$ for which $\theta_t (v^0_{1n})$ is supported in the
$(k,l)$ block subspace and is such that $Q_i \phi (v^0_{1n}) Q_j  = 0$
for all $(i,j) \neq (k,l)$ satisfying $ k \leq i$ and $ j\leq l$.
  We complete the proof
of the lemma by showing that for all $s, i, j$ the partial isometry
$\theta_s (v^0_{ij})$ is supported in the same block subspace as 
$\theta_t(v^0_{ij})$.  In particular all the summands $\theta_s$ are inner
conjugate and $\phi$ is a refinement embedding.

Note first that $ \theta_t(v^0_{11})$ is supported in the $ Q_k$ block
subspace and that $ \theta_t(v^0_{nn})$ is supported in the $ Q_l$
block subspace.  Given $ s \neq t$, there are various a priori 
possibilities for the support projections $ Q_i$ of $
\theta_s(v^0_{11})$.  The first possibility, $ i < j$, is suggested by
the following diagram:
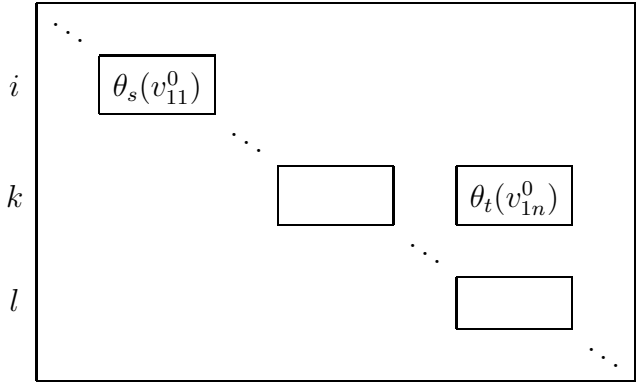
\[
\begin{array} { c | c | c | c | c | c | c | c |}
    \hhline{~-------}
   &  \multicolumn{7}{l|}{\ddots} \\
    \hhline{~~-~~~~~}
   i & \phantom{\ddots} & \theta_s(v^0_{11}) \vphantom{\ds \sum}
                & \multicolumn{5}{|c|}{} \\
    \hhline{~~-~~~~~}
   & \multicolumn{2}{c}{} &\multicolumn{1}{c}{\ddots}
    &\multicolumn{4}{c|}{} \\
    \hhline{~~~~-~-~}
   k & \multicolumn{3}{|c|}{} & \phantom{\theta_t(v^0_{1n})} 
       \vphantom{\ds \sum}
       & \phantom{\ddots} & \theta_t(v^0_{1n}) \vphantom{\ds \sum}
             & \phantom{\ddots}\\
    \hhline{~~~~-~-~}    
    & \multicolumn{4}{|c}{} &  \multicolumn{3}{l|}{\ddots}\\  
    \hhline{~~~~~~-~}
   l\vphantom{\ds \sum} & \multicolumn{5}{|c|}{} & \phantom{\theta_t(v^0_{1n})} 
       \vphantom{\ds \sum} & \\
    \hhline{~~~~~~-~}
   & \multicolumn{6}{|c}{} & \ddots  \\
    \hhline{~-------}     
\end{array}
\]
In this case we can deduce from the \orcon  nature of $\phi$
that $\theta_s(v^0_{nn}) $ is supported in a block $Q_j $ with 
$j \leq k $.  Indeed, if such a $j$ is greater than $k$, then it
follows that $\theta_s(v^0_{1n})$ and $\theta_s(v^1_{1n}) $ have
support in the $(i,j) $ block, $\theta_t(v^0_{11}) $ and
$\theta_t(v^1_{11}) $ have support in the $(k,k)$ block and so the
support of $\theta_s(v^0_{1n}) + \theta_t(v^1_{11}) $ is not of
staircase type.  Since $\theta_s + \theta_t $ is \orcon,
this is a contradiction.  

We have shown that if $i<k$ then $\theta_s + \theta_t $ is an ordered
sum.  Now suppose that $i \geq l $.  Since $Q_i $ is the support
projection for $\theta_s(v^0_{11}) $, it follows immediately (without
the need for the \orconion of $\phi $) that 
$\theta_t + \theta_s $ is an ordered sum.

Now consider the case $k < i < l $.  Then 
$\theta_s(v^1_{11}) + \theta_t(v^0_{1n}) $ has support in the $(i,i) $ 
block and the $(k,l) $ block and so is not of staircase form.  Since 
$\phi $ is \orcon, we have a contradiction once again.  Thus 
the only possible value for $i$ is $k$.

In summary, we have shown that if $s \neq t $, then
$\theta_s(v^0_{11})$ is supported in the same block $Q_k $ as
$\theta_t(v^0_{11}) $.  By the original assumption on $t$, we must
have that $\theta_s(v^0_{1n}) $ is supported in the $(k,l) $ block,
since the only alternative is that it is supported in a $(k,r) $ with 
$r > l$.  But this would imply that 
$(\theta_s + \theta_t)(v^0_{nn} + v^1_{1n}) $ is not \orcon.

For similar reasons it follows that for $1 < x < n $, the projections
$\theta_s(v^0_{xx}) $ and $\theta_t(v^0_{xx}) $ are equivalent.  For
suppose that these projections are supported in the blocks for
$Q_p $ and $Q_q $ respectively, with $p < q $.  Then
$\theta_s(v^0_{xn}) $ has support in the $(p,l) $ block while
$\theta_t(v^1_{xx}) $ has support in the $(q,q) $ block and hence
$\theta_s + \theta_t $ is not \orcon.  Similarly,
$q < p $ is not possible.

We have shown that for all $1 \leq x \leq n $ and for all $s$ and $t$
the projections $\theta_s(v^0_{xx}) $ and $\theta_t(v^0_{xx}) $ lie in 
the same block subspace and so are inner equivalent.  It follows that
$\phi$ is a refinement type embedding.
\end{proof}

\begin{cor} \label{C:struct}
Let $\phi \colon A_1 \to A_2$ be an \orcon embedding between
finite dimensional nest algebras,  where $A_1$ has  
no rank one atoms.Then $\phi$ can be written as an
ordered sum $\phi = \phi_1 + \dots + \phi_p$, where each $\phi_i$ is
either $T_2$-degenerate or of refinement type.
\end{cor}

\begin{remark}
Theorem  \ref{T:emfor} is not valid if $A_1$ has both rank one atoms 
and higher rank atoms. 
 There are examples of embeddings which are \orcon
 and order irreducible but neither $T_2$-degenerate nor
of refinement type. 
 For example, let
$\phi \colon T(2,2,1) \longrightarrow T(6,3,1)$ be given by
\[
\left(
\begin{array}{c c | c c | c}
e_1 & & & & \\
& f_1 & & & \\
\hhline{-----}
& & e_2 & & \\
& & & f_2 & x \\
\hhline{-----}
& & & & r_3
\end{array}
\right)
\longrightarrow
\left(
\begin{array}{c c c c c c | c c  c | c}
e_1 & & & & & & & & & \\
& e_1 & & & & & & & & \\
& & f_1 & & & & & & & \\
& & & f_1 & & & & & & \\
& & & & e_2 & & & & & \\
& & & & & f_2 & & & x & \\
\hhline{----------}
& & & & & & e_2 & & & \\
& & & & & & & f_2 & & x \\
& & & & & & & & r_3 & \\
\hhline{----------}
& & & & & & & & & r_3
\end{array}
\right).
\]
Only one off diagonal position has been indicated; all off diagonal
matrix units are dealt with in the same way.  Their image under $\phi$ 
matches the diagonal matrix units of the domain and the range spaces
in the order in which they appear on the diagonal.

Note that in this example, if $\phi$ is restricted to the orthogonal
complement of the rank one atom of $A_1$, then it is no longer order
irreducible. 

It is easy to construct a variety of similar counterexamples in which
the number of atoms in the codomain which intersect $\phi(1)$ is
arbitrarily large.
\end{remark}

\subsection{\Orpre embeddings.} \label{SS:more}

Donsig \cite{apd98} has extended the notion of order preservation
to the general context of digraph algebras (where order preservation
is defined with respect to a self-adjoint subalgebra of the digraph
algebra).  His definition, when 
restricted to a finite dimensional nest
algebra and taken with respect to $A \cap A^*$, agrees with the
definition below.

\begin{defn}  \label{D:sop}
Let  $A \subseteq    M_n$   be a finite dimensional nest algebra
with atomic interval projections  $q_1, q_2, \dots,  q_l$. 
\begin{enumerate}
\item 
 The set $\Pisom_{\op}( A)$ is the set of
partial isometries $v$ in $\Pisom_{\reg}(A)$ which
are \orpre  in the sense that
\[
q_i v q_j \ne 0 \implies q_s v q_t = 0, \quad \text{for }
 (s,t) \neq (i,j) \text{ with } s \leq i \text{ and } t \geq j.
\]

\item
A star extendible embedding between finite dimensional
nest algebras is {\em \orpre\/} if it 
maps \orpre  partial
isometries to \orpre partial isometries.

\item 
 A star extendible embedding $\phi \colon A_1 \to A_2$ is 
{\em \lorpre\/} if $\phi(v)$  is in 
$\Pisom_{\op}(A_2)$,
for each rank one element $v$ in $\Pisom_{\reg}(A_1)$.
\end{enumerate} 
\end{defn}

For  triangular nest algebras $T_n $, $n = 1, 2, \dots $, the
\orpre embeddings coincide with the \orcon
embeddings, since a partial isometry 
in $T_n$ is \orcon
if, and only if, it is \orpre.  In contrast, for
non-triangular algebras \orpreion is quite stringent.
For example, a star extendible embedding,
\[
\phi \colon T_2 \to T_2 \otimes M_n,
\]
is \orpre if, and only if, either
\[
\phi(e_{12}) \in D_2 \otimes M_n \quad \text{or} \quad 
  \phi(e_{12}) \in e_{12} \otimes M_n .
\]

On the other hand, if $\phi \colon T_n \to T_m $ is an \orcon
embedding and $\eta \colon B_1 \to B_2 $ is a 
\cstar algebra homomorphism between finite dimensional 
\cstar algebras, then the star extendible homomorphism,
\[
\phi \otimes \eta \colon T_n \otimes B_1 \to T_m \otimes B_2,
\]
has partial embeddings (between the summands of the algebras) which
are \lorpre.  Although the local \orpreion 
 property is not generally conserved under
composition, this is plainly the case for these tensor product maps.

We remark that there are \orcon  embeddings which are 
\lorpre but are not \orpre.  The
multiplicity two standard embedding from $T_3 \otimes M_r$ to
$T(r, r, 2r, r, r)$ has this property.

It is also possible for an embedding to be \orpre
but not \orcon. To see this consider the map 
$T_2 \otimes M_2 \rightarrow T_4 \otimes M_2$ given by
\[
\left(
  \begin{array}{c|c}
    a & b \\
    \hhline{--}
     & c \\
  \end{array}
\right) 
\longrightarrow
\left(
  \begin{array}{c|c|c|c }
    a & & b & \\
    \hhline{----}
    & a & & b \\
     \hhline{----}
     & & c & \\
      \hhline{----}
      & & & c \\
  \end{array}
\right).
\]
Finally, the embedding $T_2 \rightarrow T_2 \otimes M_3 $ given by
\[
\begin{pmatrix}
 a & b \\ 
  & c
\end{pmatrix}
\longrightarrow
\left(
  \begin{array}{c c c | c c c}
 a & b & & & & \\
   & c & & & & \\
   &   & a & b & & \\
   \hhline{------}
   & & & c & & \\
   & & & & a & b \\
   & & & & &  c
  \end{array}
\right)
\]
illustrates how easily a `nice' embedding can fail to be 
\lorpre. (For embeddings with domain $T_2$, \orpreion
and \lorpreion coincide.)

\section{Isomorphisms} \label{S:iso}

Let $\AAA = \{A_k, \alpha_k\}$, $\AAA' = \{A'_k, \alpha'_k\}$
be direct systems of digraph algebras with embeddings
$\alpha_k, \alpha'_k$ which are star extendible  
and regular.  Let $A_0, A'_0$ be the locally finite 
algebras which are the algebraic direct limits of the systems.
The limit algebras $A_0, A'_0$ are said to
be {\em star extendibly isomorphic\/} if there is an algebra
isomorphism $\Phi_0 \colon A_0 \to A'_0$ 
 which extends to a \star algebra isomorphism
between the generated C*-algebras.

\begin{definition} \label{D:sysis}
The direct systems $\AAA, \AAA'$ are said to be 
\emph{regularly isomorphic} if there exist
 regular, star extendible embeddings
$\{\phi_k, \psi_k\}$ such that the diagram
\[
\begin{diagram}
  \node{A_{n_1}} \arrow[2]{e} \arrow{se,t}{\phi_1} \node[2]{A_{n_2}}
      \arrow[2]{e} \arrow{se,t}{\phi_2} \node[2]{A_{n_3}}  
 \arrow[2]{e} \arrow{se} \\
 \node[2]{A'_{m_1}} \arrow{ne,t}{\psi_1} \arrow[2]{e}
 \node[2]{A'_{m_2}} \arrow{ne,t}{\psi_2} \arrow[2]{e} 
  \node[2]{\makebox[1 em]{$\vphantom{A_n}$}}
\end{diagram}
\]
commutes.
\end{definition}
 
A fundamental issue in the classification of nonselfadjoint limit
algebras 
is that a star extendible isomorphism 
$\Phi \colon A_0 \to A'_0$
need not be induced by   a regular
isomorphism between the given systems.
(See \cite{dp96, scp_b}.) However, we shall
see that for various \orcon  systems this is indeed the case.
This enables us to define invariants for limit algebras in terms of
regular isomorphism invariants for systems. 

In the proof of the
next lemma the block matrix support of an element $a$ in a
finite dimensional nest algebra refers to the set of block subspaces
$Q_iAQ_j$ for which $Q_iaQ_j \ne 0$.  The matrix support of $a$ refers
to a given matrix unit system $\{h_{ij}\}$ and is
the set $\{\,(i,j) \mid h_{ii}ah_{jj} \ne 0 \,\}$. 

\newcommand{\diag}[1]{A_{#1}^{\vphantom{*}} \cap A_{#1}^*}
\begin{lma} \label{L:cotmp}
    Let $\phi \colon A_1 \rightarrow A_2$, $\psi \colon A_2  
\rightarrow A_3$ be star extendible homomorphisms  between 
finite dimensional nest algebras such that the 
composition $\psi \circ  \phi$ is 
\lorcon.  If $h$ is a rank one
regular partial
isometry in $\Pisom_{\reg}(A_2)$ which lies
in the bimodule over $A_2 \cap A_2^* $ 
generated by $\phi(A_1)$,
then $\psi(h)$ is an \orcon  partial isometry.
\end{lma}

\begin{proof} Fix matrix units for $A_1$ and $A_2$ so that
$\phi$ maps matrix units of $\diag{1}$ to sums of matrix units in 
$\diag{2}$, and then choose matrix units for $A_3$ such that
$\psi$ maps matrix units in $\diag{2}$ to sums of matrix units.
Let $B$ be the bimodule over $A_2 \cap A^*_2$ generated by
$\phi(A_1)$.  Since $B$ is, in particular, a bimodule over the
diagonal matrices, $B$ is the linear span of the matrix units which it
contains. Let $e$ be a matrix unit in $B$.  Then there are atoms $P$
and $Q$ for $A_2$ such that $e = PeQ $.  If $f$ and $g$ are matrix
units satisfying $f=PfP $ and $g=QgQ $, 
then $f,g \in A_2 \cap A_2^*$.  It
follows that $feg \in B$.  But any matrix unit $e'$ for which
$e'=Pe'Q $ can be written in this form.  Thus, if $PBQ \neq \{0\} $ then 
$PBQ = PA_2Q $.  This enables us to identify the bimodule over
 $A_2 \cap A^*_2$ generated by $\phi(A_1)$: it consists of all
elements of $A_2$ which are supported in the collection of matrix
blocks which contain non-zero entries for elements of $\phi(A_1)$.

Now suppose that $h$ and $e$ are two matrix units in the same block in
$A_2$ and suppose, further, that $\psi(e)$ is \orcon.  By
the first assumption, there are matrix units $f$ and $g$ in diagonal
blocks of $A_2$ such that $h=feg $.  Since $\psi(e)$ is \orcon,
 its support has `staircase' form in $A_3$; since $\psi(f)$ 
and $\psi(g) $ lie in $A_3 \cap A_3^* $, $\psi(h)$ also has
`staircase' form; i.e.\ $\psi(h)$ is also \orcon.  This
argument works equally well if $h$ is a rank one partial isometry in
the same block as $e$ rather than a matrix unit.  (The only
modification is that $f$ and $g$ are merely rank one partial
isometries in the appropriate diagonal blocks.)  As a consequence of
these observations, we see that we can prove the lemma by proving it
in the special case in which $h$ is a subordinate of 
a partial isometry in the range, $\phi(A_1)$, of $\phi$,

If we make this additional hypothesis,
 there is a rank one regular partial
isometry $e$ in $A_1 $ such that $hh^* \phi(e) h^*h \neq 0 $.
 Because of our
matrix unit choice for $\psi$, it follows that the matrix support of
$\psi(\phi(e))$ contains the matrix support of $\psi(h)$.  However,
the block matrix support of $\psi(\phi(e))$ is of staircase type,
since $\psi \circ \phi$ is \lorcon.  It follows that
$\psi(h)$ is an \orcon partial isometry.
\end{proof}

\begin{remark}
The assumption that $\psi \circ \phi$ is \lorcon does 
not imply that $\psi$ is \lorcon, even when $\psi $
is restricted to the algebra generated by $\phi(A_1) $ and
 $A_2 \cap A_2^*$.  The two embeddings
$\phi \colon T_2 \rightarrow T(1,2,1)$ and
$\psi \colon T(1,2,1) \rightarrow T(1,1,4,1,1)$ given by
\[
\left(
\begin{array}{c | c}
  a & b \\
 \hhline{--}
  & c
\end{array}
\right)
\overset{\phi}{\longrightarrow}
\left(
\begin{array}{c | c c | c}
 a & b & & \\
 \hhline{----}
  & c & & \\
  & & a & b \\
  \hhline{----}
  & & & c
\end{array}
\right)
\]
and
\[
\left(
\begin{array}{c | c | c}
 \alpha & \beta & \gamma \\
 \hhline{---}
  & \delta & \epsilon \\
  \hhline{---}
   & & \eta
\end{array}
\right)
\overset{\psi}{\longrightarrow}
\left(
\begin{array}{c | c | c  c | c | c}
 \alpha & & & \beta & & \gamma \\
 \hhline{------}
  & \alpha & \beta & & \gamma &  \\
  \hhline{------}
  & & \delta & & \epsilon & \\
  & & &\delta & & \epsilon \\
  \hhline{------}
  & & & &  \eta & \\
 \hhline{------}
 & & & & & \eta
\end{array}
\right)
\]
provide a counter-example.
Note that $T(1,2,1)$ is the algebra generated by 
$\phi(T_2)$ and $T(1,2,1) \cap T(1,2,1)^*$.
\end{remark}

\begin{thm} \label{T:istmp}
Let ${\cal A}, {\cal A'}$ be systems of finite 
dimensional nest algebras in which all compositions of embeddings 
are \lorcon  and let $A_0, A_0'$ be their algebraic limit algebras.
Then the following are equivalent.
\begin{enumerate} \renewcommand{\theenumi}{\roman{enumi}}
\renewcommand{\labelenumi}{(\theenumi)}
\item  \label{sys} ${\cal A}$ and
 $ {\cal A'}$ are regularly isomorphic.
 \item \label{alim} $A_0 $ and $A_0'$ are 
 star extendibly isomorphic.
\end{enumerate}
\noindent
Moreover, every star extendible isomorphism $\Phi \colon A_0 \to A'_0$
is induced by a regular star extendible isomorphism
between the systems ${\cal A}, {\cal A}'$.
\end{thm}

\begin{proof}
The proof that (i) implies (ii) is immediate.  For the converse,
let $ \Phi \colon A_0 \to A'_0$ be a star extendible isomorphism.
Since each algebra in the systems $ \cal A$ and $ \cal A'$ is finitely 
generated, there is a commuting diagram isomorphism of the
systems with crossover maps $\phi_k $, $\psi_k $:
\[
\begin{diagram}
  \node{A_{n_1}} \arrow[2]{e,t}{i_1} \arrow{se,t}{\phi_1}
 \node[2]{A_{n_2}}
      \arrow[2]{e,t}{i_2} \arrow{se,t}{\phi_2} \node[2]{A_{n_3}}  
 \arrow[2]{e,t}{i_3} \arrow{se} \\
 \node[2]{A'_{m_1}} \arrow{ne,t}{\psi_1} \arrow[2]{e,t}{i'_1}
 \node[2]{A'_{m_2}} \arrow{ne,t}{\psi_2} \arrow[2]{e,t}{i'_2} 
  \node[2]{\makebox[1 em]{$\vphantom{A_n}$}}
\end{diagram}
\]

Since $i_1 = \psi_1 \circ \phi_1 $, the range of $i_1 $ is contained
in $\psi(A'_{m_1}) $.  By hypothesis, $i_1$ is \lorcon
 and so locally regular and hence
regular.  If $x \in A_{n_1}$ is a rank one
regular partial isometry, then $i_1(x) $ is a regular, \orcon,
 partial isometry in $\psi(A'_{m_1}) $.  Consequently,
$i_1(x) $ is a sum, $i_1(x) = y_1 + \dots + y_s$, of regular rank one
partial isometries in $\psi(A'_{m_1})$.  Since 
$i'_1 = \phi_2 \circ \psi_1$ is \lorcon,
 each partial isometry $\phi_2(y_j)$ is \orcon,
by Lemma~\ref{L:cotmp}, and hence
 regular in $A'_{m_2}$.  Therefore
$\phi_2 \circ i_1 (x) $ is regular and so 
$\phi_2 \circ i_1 \colon A_{n_1} \to A'_{m_2} $ is locally regular.
By Theorem~\ref{T:reg}, 
$\phi_2 \circ i_1 = \phi_2 \circ \psi_1 \circ \phi_1$ is a regular
embedding. 

The same argument shows that 
$\psi_3 \circ \phi_3 \circ \psi_2 = \psi_3 \circ i'_2$ is regular;
continuing in this fashion, we may replace the initial commuting
diagram isomorphism with one whose crossover maps are all regular.
This shows that (ii) implies   (i).
\end{proof}

The last theorem is important in that it shows that tentative invariants
defined for algebraic limit algebras in terms of specific regular
presentations do not in fact depend on the presentation and are
thus genuine invariants for star extendible isomorphism.
In the next section we shall need to know that certain scales defined
in terms of \orpre or  \orcon systems are invariants. This will follow
from the following refinement of Theorem~\ref{T:istmp}
in the order preserving or \orcon cases.

\begin{thm} \label{T:autooc}
Let ${\cal A}, {\cal A'}$ be systems of finite 
dimensional nest algebras in which all  embeddings 
are \orpre (resp. \orcon)  and let $A_0, A_0'$ be 
their algebraic limit algebras.
Then the following are equivalent.
\begin{enumerate} \renewcommand{\theenumi}{\roman{enumi}}
\renewcommand{\labelenumi}{(\theenumi)}
\item  \label{sys2} The systems ${\cal A}$ and
 $ {\cal A'}$ are regularly isomorphic by a commuting diagram 
isomorphism in which all the crossover maps are 
\orpre (resp. \orcon). 
 \item \label{alim2} $A_0 $ and $A_0'$ are 
 star extendibly isomorphic.
\end{enumerate}
\noindent
Moreover, every star extendible isomorphism $\Phi \colon A_0 \to A'_0$
is induced by a system isomorphism as in \textup{(\ref{sys2})}.
\end{thm}

\begin{proof}
The proof of this theorem follows easily from
 Lemma~\ref{L:autooc} below.  If $A_0$ and $A_0'$ are star extendibly
 isomorphic, then Theorem~\ref{T:istmp} gives a regular system
 isomorphism in which any adjacent pair of crossover maps have a
 composition which is \orpre or \orcon, as appropriate.  From the
 lemma, any consecutive triple of crossover maps has composition which
 is \orpre or \orcon; a sequence of such triple compositions yields a
 system isomorphism with \orpre or \orcon crossover maps.
\end{proof}

\begin{lemma} \label{L:autooc}
Let
\[
A_1 \overset{\phi}{\longrightarrow} A_2
    \overset{\psi}{\longrightarrow} A_3
    \overset{\eta}{\longrightarrow} A_4
\]
be regular star extendible homomorphsims between finite dimensional
nest algebras such that the compositions $\psi \phi$
and $\eta \psi$ are \orpre (resp. \orcon).
Then the triple composition  $\eta \psi \phi$
is \orpre (resp. \orcon).
\end{lemma}

Before giving a proof of the Lemma, we discuss some helpful
preliminary matters.
When $A_1$ is a finite dimensional nest, we number the block rows and
block columns of $A_1$ in the natural way.  If $v \in A_1$ is
supported in a single block, $\row(v)$ and $\col(v)$ denote the row
and column in which the support block is located.  Suppose that
$\phi \colon A_1 \to A_2$ is a multiplicity one 
(star extendible, regular) embedding of $A_1$
into another finite dimensional nest algebra.  If $v$ is supported in
a single block in $A_1$, then $\phi(v)$ is supported in a single block
in $A_2$.  If $v$ and $w$ are each supported in a single block, then
\begin{align*}
\row(v) = \row(w) &\Longrightarrow \row(\phi(v)) = \row(\phi(w)), \\
\col(v) = \col(w) &\Longrightarrow \col(\phi(v)) = \col(\phi(w)), \\
\row(v) < \row(w) &\Longrightarrow \row(\phi(v)) \leq \row(\phi(w)),\\
\col(v) < \col(w) &\Longrightarrow \col(\phi(v)) \leq \col(\phi(w)).
\end{align*}  
In particular, if $v$ and $w$ are supported in the same block, then
so are $\phi(v)$ and $\phi(w)$.
Thus, $\phi$ induces a map $\tph$ from the blocks of $A_1$ to the
blocks of $A_2$.  This map need not be injective, but it does respect
the block structure:  if $\row(X) < \row(Y)$, then
$\row(\tph(X)) \leq \row(\tph(Y))$; 
similarly, if $\col(X) < \col(Y)$, then
$\col(\tph(X)) \leq \col(\tph(Y))$.  (Here, $X$ and $Y$ denote blocks
in $A_1$, not elements of the nest algebra.)
In particular, if $X$ and $Y$ are blocks in
the same  row (or  column), then $\tph(X)$ and $\tph(Y)$ are
also blocks in the same  row (or column).

Let $v$ and $w$ be partial isometries each of which is supported in a
single block and assume that $v+w$ is a partial isometry which is not
\orpre. Then the block support for one of the partial
isometries, say $v$, is located to the ``southwest'' of the block
support for the other partial isometry.  More precisely,
\[
\row(w) \leq \row(v) \quad \text{and} \quad \col(v) \leq \col w,
\]
with at least one inequality being strict.  This can be indicated by
the diagrams:
\[
\nofusion{\phantom{w}}{}{v}{w}\,,\quad 
\nofusion{w}{}{v}{\phantom{w}}\,,\quad
\text{or}\quad \nofusion{}{w}{v}{\phantom{w}}\,.
\] 

Now suppose that $\phi(v+w) = \phi(v) + \phi(w)$ is \orpre.
 (The assumption that $\phi$ has multiplicity one is still
in force.) It follows from the relations above that the block support
for $\phi(v)$ is to the ``southwest,'' i.e., that
\[
\row \phi(w) \leq \row \phi(v) \quad \text{and} \quad
\col \phi(v) \leq \col \phi(w).
\]
But $\phi(v) + \phi(w)$ is \orpre, so in fact
$\row \phi(v) = \row \phi(w)$ and $\col \phi(v) = \col \phi(w)$;
i.e., $\phi(v)$ and $\phi(w)$ have the same block support.

Now assume that $\phi \colon A_1 \to A_2$ is an embedding with
arbitrary multiplicity.  Suppose that $\phi$ is not \orpre.
Then one (or both) of two situations must occur.  The first is that
there is a rank one partial isometry $v \in A_1$ and a pair of
multiplicity one summands $\phi_1$ and $\phi_2$ of $\phi$ such that
$\phi_1(v) + \phi_2(v)$ is not \orpre.  The other
alternative is that there are two rank one partial isometries $v$ and
$w$ in $A_1$ with support in distinct blocks such that $v+w$ is 
\orpre and there are two multiplicity one summands 
 $\phi_1$ and $\phi_2$ of $\phi$ such that $\phi_1(v) + \phi_2(w)$ is
 not \orpre.

If, on the other hand, $\phi$ is \orpre, then any (partial)
sum of the multiplicity one summands for $\phi$ is also \orpre.
Likewise, if  $\phi$ is \orcon, then any (partial)
sum of the multiplicity one summands for $\phi$ is also \orcon.

Note also that a multiplicity one summand of a composition of two
embeddings is the composition of multiplicity one summands of each
factor. 

\begin{proof}[Proof of Lemma~\ref{L:autooc}] 
The proof of the Lemma in the \orpre context differs considerably 
from the proof in the \orcon context, so we present the two arguments
separately.   We start with the
the \orpre context, where the argument is simpler.

First, suppose that $v$ is a rank one partial isometry in $A_1$ such
that $\eta \psi \phi (v)$ is not \orpre.
Then there are
multiplicty one summands $\phi_1$ and $\phi_2$ of $\phi$, 
$\psi_1$ and $\psi_2$ of $\psi$, and $\eta_1$ and $\eta_2$ of $\eta$
such that $\eta_1 \psi_1 \phi_1 (v) + \eta_2 \psi_2 \phi_2 (v)$ is not
\orpre.

Observe that $\phi_1(v) + \phi_2(v)$ is a partial isometry which
is not \orpre.  
It is a partial isometry since $\phi_1 + \phi_2$ is an embedding; 
if it were \orpre, then its image under 
$\eta_1 \psi_1 + \eta_2 \psi_2$ would be \orpre. 
 But $\eta_1 \psi_1 \phi_1 (v) + \eta_2 \psi_2 \phi_2 (v)$ 
is a subordinate of this image and is not \orpre. 

Since $\psi_2 \phi_1 + \psi_2 \phi_2$ (a sum of multiplicity one
summands of $\psi \phi$) is \orpre, it follows that
$\psi_2 \phi_1 (v) + \psi_2 \phi_2 (v)$ 
is \orpre.  Consequently,
$\psi_2 \phi_1 (v)$ and $\psi_2 \phi_2 (v)$ are in the same block.
This, in turn, implies that $\eta_2 \psi_2 \phi_1 (v)$ and
$\eta_2 \psi_2 \phi_2 (v)$ are in the same block.

This means that
\[
\eta_1 \psi_1 \phi_1 (v) + \eta_2 \psi_2 \phi_1 (v)
\]
and
\[
\eta_1 \psi_1 \phi_1 (v) + \eta_2 \psi_2 \phi_2 (v)
\]
have the same block structure -- either both are \orpre or
both are not \orpre.  But the second one is not \orpre
 by assumption while the first one is the image of the rank
one partial isometry $\phi_1(v)$ under the \orpre embedding
$\eta_1 \psi_1 + \eta_2 \psi_2$, a contradiction.

The other possibility which we need to consider is that there are
rank one partial isometries $v$ and $w$ in $A_1$ such
$v+w$ is an \orpre partial isometry but
 $\eta \psi \phi (v+w)$ is not \orpre.
Again, there are
multiplicty one summands $\phi_1$ and $\phi_2$ of $\phi$, 
$\psi_1$ and $\psi_2$ of $\psi$, and $\eta_1$ and $\eta_2$ of $\eta$
such that $\eta_1 \psi_1 \phi_1 (v) + \eta_2 \psi_2 \phi_2 (w)$ is not
\orpre.

Observe that $\phi_1(v) + \phi_2(w)$ is a partial isometry which
is not \orpre.  
It is a partial isometry since it is a subordinate of
$\phi(v+w)$;
if it were \orpre, then its image under 
$\eta_1 \psi_1 + \eta_2 \psi_2$ would be \orpre. 
 But $\eta_1 \psi_1 \phi_1 (v) + \eta_2 \psi_2 \phi_2 (w)$ 
is a subordinate of this image and is not \orpre. 

Since $\psi_2 \phi_1 + \psi_2 \phi_2$ 
 is \orpre, it follows that
$\psi_2 \phi_1 (v) + \psi_2 \phi_2 (w)$
 is \orpre.  Consequently,
$\psi_2 \phi_1 (v)$ and $\psi_2 \phi_2 (w)$ are in the same block,
whence $\eta_2 \psi_2 \phi_1 (v)$ and
$\eta_2 \psi_2 \phi_2 (v)$ are in the same block.

This means that
\[
\eta_1 \psi_1 \phi_1 (v) + \eta_2 \psi_2 \phi_1 (w)
\]
and
\[
\eta_1 \psi_1 \phi_1 (v) + \eta_2 \psi_2 \phi_2 (w)
\]
have the same block structure  But 
$\eta_1 \psi_1 \phi_1 (v) + \eta_2 \psi_2 \phi_2 (w)$
  is not \orpre by assumption while 
$\eta_1 \psi_1 \phi_1 (v) + \eta_2 \psi_2 \phi_1 (w)$
is a subordinate of the image of the \orpre
partial isometry $\phi_1(v+w)$ under the \orpre embedding
$\eta_1 \psi_1 + \eta_2 \psi_2$, a contradiction.
($\phi_1(v+w)$ is \orpre since any multiplicity one
embedding is \orpre and $v+w$ is an \orpre partial
isometry.)

Having completed the \orpre context, we move on to the \orcon
context. 
Let $v \in A_1$ be an \orcon partial isometry such that the
block support of $v$ has staircase form.
Let $v' = v_a + v_b$
be a rank two subordinate of $v$ which is supported in two of these
block subspaces, denoted $a$ and $b$. 
The orientation of these blocks in the block decomposition of
$A_1$ can
be indicated diagramatically as one of six types, namely the triple,
\[
\nofusion{a}{\phantom{b}}{}{b}\,, \qquad
\nofusion{a}{b}{}{}{}\,, \qquad
\nofusion{}{b}{\phantom{b}}{a}\,, \qquad
  \]
together with the corresponding triple with the letters reversed.

Suppose the $\phi(v')$ is not \orcon. Then there exist two
blocks of $A_2$,  $X$ and $Y$ say, which have nonstaircase orientation
\[
 \nofusion{}{}{X}{Y} 
 \]
and are such that at least one of four possibilities occurs:

(i) $\phi(v_a)$ has support meeting $X$ and $Y$

(ii) $\phi(v_b)$ has support meeting $X$ and $Y$

(iii)  $\phi(v_a)$ and $\phi(v_b)$ have  support meeting $X$ and $Y$
respectively

(iv)  $\phi(v_a)$ and $\phi(v_b)$ have  support meeting $Y$ and $X$
respectively.

Consider the companion blocks $U$, $V$ for $X$, $Y$
indicated by  the diagram
\[
\nofusion{U}{V}{X}{Y}\,.
\]
Since the composition $\psi \phi$ is \orcon,
 if $\psi_1$ is a multiplicity one
summand of $\psi$, 
the composition $\psi_1 \phi$ is \orcon.
In particular such a map $\psi_1$ must `correct' the $X$, $Y$ support of
$\phi(v)$
(indicated in one of the four possibilities (i) to (iv)) by
a `fusion' of blocks, as indicated in each of the three diagrams
\begin{align*}
 &\nofusion{U}{V}{X}{Y} \longrightarrow 
\hfusion{\tps(U)}{\tps(V)}{\tps(X)}{\tps(Y)}\,, \\
&\nofusion{U}{V}{X}{Y} \longrightarrow 
\vfusion{\tps(U)}{\tps(V)}{\tps(X)}{\tps(Y)}\,, \\
&\nofusion{U}{V}{X}{Y} \longrightarrow 
\fullfusion{\tps(U)}{\tps(V)}{\tps(X)}{\tps(Y)}\,.
\end{align*}
This means that, in the first case for example, the image under
$\psi_1$ of an element in $A_2$ with support in the blocks $U$, $Y$
or in the blocks $X$, $V$ has support in
a single block of $A_3$.  (The appearance of $\tps(U)$ and $\tps(Y)$,
for example, in the same block in the diagram indicates that
$\tps(U) = \tps(Y)$.)
 These types of fusion may vary among
the various summands 
$\psi_1$ of $\psi$.
The comments preceding the statement of the lemma preclude a
correction of the $X$, $Y$ support of $\phi(v)$ of the form
\[
\nofusion{}{}{X}{Y} \longrightarrow 
\nofusion{\tps(X)}{\tps(Y)}{}{}
\]
or the similar correction with $\tps(X)$ and $\tps(Y)$
interchanged.

Let $w, u \in A_2$  be rank two partial isometries with support in the
blocks $X$, $Y$ and $U$, $V$ respectively. We now show that 
$\eta \psi(w)$ is \orcon.
The idea for this argument is that the fusion of blocks `binds' the
supports of $\eta \psi(w)$
and $\eta \psi(u)$ and the latter is \orcon by the
hypotheses, since $u$ is \orcon. More precisely let $w = w_X
+ w_Y$ be the rank one decomposition and suppose that 
$\eta \psi(w)$ is not \orcon.
Then there exist two blocks $W, Z$ in $A_4$ with nonstaircase form
\[
 \nofusion{}{}{W}{Z}  
\]
and there exist multiplicity one summands $ \psi_1$, $\psi_2$ of
$\psi$ and $\eta_1$, $\eta_2$ of $\eta$ such that at least one of the
following possibilities occurs.

(a) $\eta_1 \psi_1(w_X)$ and  $\eta_2 \psi_2(w_X)$
meet $W$ and $Z$ respectively.

(b)  $\eta_1 \psi_1(w_Y)$ and $ \eta_2 \psi_2(w_Y)$
meet $W$ and $Z$ respectively.

(c)  $\eta_1 \psi_1(w_X)$ and $ \eta_2 \psi_2(w_Y)$
meet $W$ and $Z$ respectively.

(d)  $\eta_1 \psi_1(w_X) $ and $ \eta_2 \psi_2(w_Y)$
meet $Z$ and $W$ respectively.

Let $u = u_U + u_V$ be the rank one decomposition of $u$.
If the first possibility (a) occurs, then noting that $\psi_1$ fuses
$X$ and $U$, or $X$ and $V$, and that $\psi_2$ also fuses
$X$ and $U$, or $X$ and $V$, it follows 
that $\eta_1 \psi_1(u_U)$
(or  $\eta_1 \psi_1(u_V)$) meets $W$ and similarly that
$\eta_2 \psi_2(u_U)$
(or  $\eta_2 \psi_2(u_V)$) meets $Z$.
For each of these four alternatives 
$\eta_1 \psi_1(u) + \eta_2 \psi_2(u)$
is not \orcon, and this
 contradicts the fact that $\eta \psi$
 is \orcon. The other three possibilities 
(b), (c), (d) also lead to
contradictions in the same manner.

It has been shown then that $\eta \psi(w)$
 is \orcon or, more intuitively, that the map
$\eta \psi$ `corrects' the $X$, $Y$ block structure of $\phi(v')$ and
$\phi(v)$. We now wish to deduce that since every such block pair is
corrected by $\eta \psi$ then in fact $\eta \psi \phi(v)$ is \orcon.
 This will complete the proof.

To see this, suppose that $\eta \psi \phi(v)$
is not \orcon and that it has support in two blocks $P$, $Q$
of $A_4$ which are not in staircase form.
Then there are multiplicity one summands 
$\eta_1 \psi_1 \phi_1$ and $\eta_2 \psi_2 \phi_2$
of $\eta \psi \phi$ (where $\eta_1$, $\eta_2$ are multiplicity
one summands of $\eta$; $\psi_1$, $\psi_2$ are multiplicity
one summands of $\psi$; and $\phi_1$, $\phi_2$ are
multiplicity one summands of $\phi$)
 such that $\eta_1 \psi_1 \phi_1(v)$ meets $P$ 
and $\eta_2 \psi_2 \phi_2(v)$ meets $Q$. 
Thus there are rank one subordinates $v_1, v_2$ of $v$ such that
$\eta_1 \psi_1 \phi_1(v_1)$ meets $P$ and
 $\eta_2 \psi_2 \phi_2(v_2)$ meets $Q$. 
It follows that the partial isometry $w = \phi_1(v_1) + \phi_2(v_2)$
is not \orcon (since the map
 $\eta \psi$ is \orcon).  However the argument 
above shows that $\eta \psi(w)$ is \orcon, which
is the desired  contradiction.
  \end{proof}

\begin{remark}
The proofs of Theorems~\ref{T:istmp} and~\ref{T:autooc} carry 
over readily to the case of
systems of direct sums of finite dimensional nest algebras for which
the embeddings have \orpre or \orcon partial embeddings between nest
algebra summands.
\end{remark}

\begin{remark}
In the terminology of Power \cite{scp_b} 
Theorem~\ref{T:autooc} establishes that the families of order
preserving and order conserving embeddings are functorial families. It
follows from Haworth and Power \cite{hawp_a} that standard
AF diagonal masas are unique up to approximately inner automorphism
in the corresponding algebraic limit algebras.
In particular the spectrum is a well defined and a complete invariant
for these algebras.
\end{remark}

\section{Invariants and classifications} \label{S:invar}

Let $A$ be a digraph algebra whose block diagonal  
subalgebra $A \cap A^*$ has minimal central projections
$q_1, q_2, \dots , q_n $.  The reduced graph $H$ for $A$ is the graph
with vertices labelled $1, \dots, n $ and edges $(i,j)$ for which
$q_i A q_j $ is non-zero.  The abelian group $\bbz^{\omega(A)} $,
where $\omega(A) $ is the number of edges of $H$, may be defined
intrinsically as the free abelian group $\GG(A) $ whose generators are 
the equivalence classes $[v]$ of rank one elements $v$ of
$\Pisom_{\reg}(A)$, where the equivalence relation is the following:
$v \sim w$ if, and only if, there are unitaries $u_1$ and $u_2$ in
$A$ such that $v = u_1wu_2$.
  (Recall that $v \in \Pisom_{\reg}(A)$ 
requires that $v^*v $ and 
$vv^* $ belong to $A$.)  If $A$ is a finite dimensional nest algebra and 
$q_1, q_2, \dots, q_n $ are the atoms of $A$, in order, then we write
$T_n(\bbz) $ to denote the abelian
group  $\GG(A) $.  Note that the group homomorphisms
$\pi_f \colon \GG(A) \to K_0(A)$ and $\pi_i \colon \GG(A) \to K_0(A)$
determined by the correspondences $[v] \mapsto [vv^*] $ and
 $[v] \mapsto [v^*v] $ respectively correspond to the row sum
 homomorphism and the column sum homomorphism from $T_n(\bbz) $  to 
$\bbz^n $.  Observe also that $K_0(A)$ can be identified in a natural
way with a subgroup of $\GG(A)$.  For a full treatment of this
invariant, see \cite{scp}.

\begin{definition} \label{D:ddgr}
The \emph{dimension distribution group} of the digraph algebra
 $A$ is the group $\GG(A) $
described above.  If $\cal A = \{ A_k, \alpha_k \} $ is a regular star 
extendible system of digraph algebras, the \emph{dimension distribution
group}, $\GG({\cal A})$, of $\cal A$ is
defined to be  the direct limit
$\varinjlim(\GG(A_k), \GG(\alpha_k)) $ determined by the naturally
induced embeddings.
\end{definition}

Plainly, only the local regularity of $\alpha_k $ 
is  necessary for the existence of the induced group homomorphisms 
$\GG(\alpha_k)$.  In view of the commuting diagram
\[
\begin{diagram}
   \node{\GG(A_k)} \arrow{s,l}{\pi} \arrow{e,t}{\GG(\alpha_k)}
     \node{\GG(A_{k+1})} \arrow{s,l}{\pi} \\
   \node{K_0(A_k)} \arrow{e,t}{K_0 (\alpha_k)} \node{K_0(A_{k+1})} 
\end{diagram}
\]
with $\pi = \pi_i $ or $\pi = \pi_f $, the invariant $\GG(\cal A)$  
comes equipped with group homomorphisms
 $\pi_i \colon \GG(\cal A) \to K_0(\cal A) $
and $\pi_f \colon \GG(\cal A) \to K_0(\cal A) $ which are, 
in fact, scaled group 
homomorphisms with respect to the $K_0 $ scale $\Sigma_0(\cal A) $ and 
the dimension group scale $\Sigma(\cal A) \subseteq \GG(\cal A) $,
which is determined by the regular partial isometries of
the system  $\cal A $.
The natural full dimension distribution group invariant is therefore
the quadruple 
\[
\Inv(\cal A) = \left( (\GG(\cal A), \Sigma(\cal A)),
   (K_0 \cal A, \Sigma_0 \cal A), \pi_f, \pi_i \right),
\]
which is an invariant for regular star extendible isomorphism.

For self-adjoint systems, $\Inv(\cal A) $ reduces to 
$(K_0 \cal A, \Sigma_0 \cal A) $. At the other extreme of
triangular algebra systems (with regular embeddings), $\Inv(\cal A)$ 
can be identified with a dual form of the fundamental binary relation
(or spectrum) invariant. See \cite{scp}. 
In the discussion below we extend the scope
for $\Inv$ as a classifying invariant.

In \cite{scp} it was shown that partly self-adjoint $T_r$-algebra
systems are classified by $\Inv $, together with an extra
 matrix unit scale, for $r = 2,3 $ but not for $r \geq 4 $.  Also, it
 was shown how one can introduce appropriate Grothendieck group
 invariants which capture the extra variety of regular embeddings and
 which serve as complete classifying invariants (at least for regular
 isomorphisms of systems).  However, it is natural to enquire to what
 extent $\Inv(-)$ already serves as a basis for complete classification
 for systems of restricted type, such as \orcon embeddings.
 This is particularly so if it is known that the inner conjugacy class of
 an admissible embedding is determined by the induced map on
 $K_0(-) $ or  on $\Inv(-) $.

\subsection{Existence.} \label{SS:exis}

Recall first the usual scheme for proving that an invariant 
such as $\Inv(-) $ is, in fact, 
a complete invariant for star extendible
isomorphism. Start with an isomorphism
$\Phi \colon \Inv(\cal A) \to \Inv(\cal A') $; this gives rise to a 
commuting diagram with crossover maps
\[
\Inv(A_{n_1}) \overset{\gamma_1}{\longrightarrow}
  \Inv(A'_{m_1}) \overset{\delta_1}{\longrightarrow}
  \Inv(A_{n_2}) \longrightarrow \dots .
\]
Lift $\gamma_1 $ to an admissible embedding 
$\phi_1 \colon A_{n_1} \to A'_{m_1}$, perhaps after increasing $m_1$.
This step, the existence step, is usually the crux of the proof and
some type of scale preservation is usually needed to enable the
lifting.  Similarly, lift $\delta_1 $ to an admissible embedding
$\psi_1 $.  Provided that admissible embeddings are closed under
composition or that, for some other reason, the composition is
admissible, appeal  next to uniqueness to conclude that 
$\psi_1 \circ \phi_1 $ is conjugate to the given embedding
$A_{n_1} \to A_{n_2} $.  Adjust $\psi_1 $ and  obtain the commuting
triangle for the first stage of the desired infinite commuting
diagram. 

We now  obtain some existence results in the context of \orconion.

For a finite dimensional nest algebra $A$, write 
$\Sigma_{\oc}(A) $ for the \orcon elements 
in the scale $\Sigma(A)$
of $\GG(A)$.  If $g$ lies in $\Sigma_{\oc}(A) $ then, under the
identification of $\GG(A)$ with $T_r(\bbz)$, $g$ is a non-negative
integral matrix and has support in staircase form. 
Conversely, if $g$ has staircase form and $\pi_i(g) \in \Sigma_0(A)$
and $\pi_f(g) \in \Sigma_0(A)$ then $g \in \Sigma_{\oc} (A)$.
 In a similar way,
define $\Sigma_{\op}(A) $, the subset of \orpre
elements. 
 
We are interested in lifting a $\pi$-respecting group homomorphism 
$\gamma \colon \GG(A_1) \to \GG(A_2) $ to a star extendible regular
embedding.  The possible \orconion and \orpreion
properties which we might
require of $\gamma$ are that it be \lorcon, \lorpre,
\orcon,  or \orpre.
(\Orconion  or \orpreion  for
$\gamma$ means, of course, that $\gamma$ maps $\Sigma_{\oc}$ into
$\Sigma_{\oc}$ or $\Sigma_{\op}$ into $\Sigma_{\op}$, respectively.
The local versions mean that $\gamma$ maps elements of $\Sigma_{\oc}$
which have only one non-zero entry into $\Sigma_{\oc}$ or
$\Sigma_{\op}$, as appropriate.)
Of course, if $\gamma$ can be lifted to a regular
embedding $\phi$ with $\GG(\phi) = \gamma $, then $\phi $ will
necessarily have the same \orconion and \orpreion
 properties as $\gamma$.
Since an embedding of any of these types is necessarily \lorcon,
it follows from Theorem~\ref{T:lift} that \lorpre
 embeddings, \orpre embeddings, and \orcon
 embeddings all have liftings of the corresponding 
 type.

\begin{thm} \label{T:lift}
Let $A_1$ and $A_2$ be finite dimensional nest algebras.
Let $\gamma \colon \GG(A_1) \to \GG(A_2) $ be a \lorcon
 homomorphism which is $\pi$-respecting and for which the
restriction of $\gamma$ to $K_0(A_1) $ is a scaled group
homomorphism.  Then there is a \lorcon star
extendible embedding $\phi \colon A_1 \to A_2 $ for which
$\GG(\phi) = \gamma$.
\end{thm}

While the proof of the theorem is somewhat technical, 
 the existence of liftings of \lorpre
embeddings is quite elementary for the following reason.  With the
obvious notation, note that one can lift the images
$\gamma([e_{12}]), \gamma([e_{23}]), \dots, \gamma([e_{r-1,r}]) $ to
\orpre regular partial isometries
$v_1, v_2, \dots, v_{r-1} $ in $A_2$ with matching initial and final
projections 
\[
v_1^*v_1 = v_2 v^*_2, v^*_2 v_2 = v_3 v^*_3, \dots, 
  v^*_{r-2} v_{r-2} = v_{r-1} v^*_{r-1}.
\]
(This uses the fact that $\gamma$ preserves the scale and respects the
group homomorphisms $\pi_f$ and $\pi_i$.)
Since the \orpre partial isometries form a
semigroupoid (admissible products of \orpre partial isometries are
\orpre) we can define $\phi$ 
as the unique star extendible embedding
for which $\phi(e_{i, i+1}) = v_i $, $1 \leq i \leq r-1$.
Now $\gamma$ and $G(\phi)$ are \lorpre
homomorphisms which are $\pi$-respecting and agree on 
$[e_{11}], \dots , [e_{rr}]$.
Since $\pi_i$ and $\pi_f$ determine an \orpre
 element it follows that
$\gamma = G(\phi)$.

This approach is not available for \lorcon
embeddings (since a product of two \orcon partial
isometries can be a partial isometry which is not \orcon).
 However,  Theorem~\ref{T:lift} follows quickly from the
following lemma.  To clarify the formulation of the 
lemma we remark that it is used to lift $\gamma([e_{12}])$ 
to a partial isometry which has \orcon products with a
priori liftings of $\gamma([e_{ij}])$, for $i \geq 2 $.
This in turn enables an immediate  proof of Theorem~\ref{T:lift}
 by induction.

\begin{lemma} \label{L:lift}
Let $X = (x_{ij})$, $Y^k = (y^k_{ij})$, and $Z^k = (z^k_{ij})$, 
$1 \leq k \leq t $, be matrices with non-negative integral entries and 
with \orcon staircase form.  Suppose that $X$ is
$p \times q $ , $Y^k$ is $q \times r_k $, $Z^k$ is
$p \times r_k$ and that $\pi_f(X) = \pi_f(Z^k) $, 
$\pi_i(X) = \pi_f(Y^k) $, and $\pi_i(Y^k) = \pi_i(Z^k) $, 
for $1 \leq k \leq t $.  Suppose moreover that 
$v_1, \dots, v_t $ are block matrices which are regular partial
isometries, with common final projections, for which each $v_k$ has
rank distribution $Y^k $.  Then there is a regular partial isometry
$u$, with rank distribution $X$ and initial projection equal to
$v_1 v_1^* $, for which the products $uv_k$ are regular partial
isometries with rank distribution $Z^K$, for $1 \leq k \leq t $.
\end{lemma} 

\begin{proof}
It will be enough to show that the assertion of the lemma follows if
it is assumed that the lemma is true in the case of matrices of sizes
$p' \times q' $, $q' \times r'_k $ and $p' \times r'_k $, where
$p' \leq p $, $q' \leq q $, $r'_k \leq r_k $, $1 \leq k \leq t $, and
at least one of these inequalities is strict.

We may assume that all the entries $x_{11}$, $y^k_{11} $ and
 $z^k_{11}$, $1 \leq k \leq t $ are non-zero.  The reason for this is
 that if a staircase form matrix has its first entry equal to zero,
 then either the first row or the first column of the matrix is zero;
 the assertion then follows from the induction hypothesis.

Let $d = \min\{x_{11}, y^k_{11}, z^k_{11} \mid 1 \leq k \leq t  \}$.
Let $e$ be a subprojection of 
$v_1 v^*_1 ( =  v_i v^*_i  \text{ for all } i)$ which has rank $d$; let 
$v'_k = ev_k $; and let $u'$ be a regular partial isometry of rank $d$ 
with initial projection $e$ and with $p \times q $ rank distribution
matrix of the form:
\[
\begin{bmatrix}
 d & 0 &  \dots&  0 \\
 0 & 0 &  \dots&  0 \\
 \vdots & \vdots & & \vdots \\
 0 & 0 &  \dots&  0 \\
\end{bmatrix}
\]
Set $z'_k = u' v'_i$  ($= u'v_i$), so that $z'_k $ has a similar rank
distribution matrix, of size $p \times r_k $.  Consider now the
matrices $\hat{X} $, $\hat{Y}_k $, and $\hat{Z_k} $, which are
obtained from $X$, $Y_k$, and $Z_k$ by subtracting $d$ form the first
entry, together with the partial isometries $\hat{v}_k = v_k - v'_k $
and notice that these satisfy the hypotheses of the lemma.  By the
induction hypothesis, we may assume that there is a lifting $\hat{u}$
for $\hat{X}$ with $\hat{u}^* \hat{u} = \hat{v}_k \hat{v}^*_k$, for
all $k$, such that $\hat{u} \hat{v}_k$ has rank distribution
$\hat{Z}_k $. Now  $u = u' + \hat{u} $ is a partial isometry with the
properties needed to prove the lemma.
\end{proof}

\subsection{Classifications} \label{SS:class}

We are now in a position to obtain a variety of classification
results. Let $\sF_{\loc}$  be the family of embeddings between finite
dimensional nest algebras which are \lorcon and let
$\Sys (\sF_{\loc})$ be the family of systems $\sA$ for which the given
embeddings and their compositions are \lorcon.
Similarly define the families
$\Sys (\sF_{\oc})$,
$\Sys (\sF_{\lop})$ and
$\Sys (\sF_{\op})$ for \orcon, \lorpre and
\orpre embeddings.

Recall that an 
isomorphism $\gamma \colon \Inv(\AAA) \longrightarrow \Inv(\AAA')$
is an isomorphism of $\GG(\AAA)$ onto $\GG(\AAA')$ which respects 
$K_0$, the scales, $\pi_i$, and $\pi_f$.  More specifically, 
$\gamma \colon \GG(\AAA) \longrightarrow \GG(\AAA')$ satisfies:
\begin{align*}
\gamma(\Sigma(\AAA)) &= \Sigma(\AAA'), \\
\gamma(K_0(\AAA)) &= K_0(\AAA'), \\
\gamma(\Sigma_0(\AAA)) &= \Sigma_0(\AAA'), \\
\gamma \circ \pi_i &= \pi_i \circ \gamma, \\ 
\gamma \circ \pi_f &= \pi_f \circ \gamma. \\ 
\end{align*}

In addition to the scale $\Sigma(\AAA)$ in $G(\AAA)$ 
 we define the following
subscales.

\begin{definition}
The scale $\Sigma_{\oc}(\AAA)$ (resp. $\Sigma_{\op}(\AAA)$) consists of the
images of the classes $[v]$ in $G(A_k)$ for which $[i(v)]$ belongs to
$\Sigma_{\oc}(A_n)$ (resp. $\Sigma_{\op}(A_n)$)  for all the system maps
$i \colon  A_k \to A_n$, for $n > k$.
\end{definition}

We now define \orconion and \lorconion
 for scale preserving homomorphisms between the dimension
distribution groups of systems.  This is, of course,
based on the corresponding definitions in the context of dimension
distribution groups of digraph algebras (see the paragraph preceding
Theorem~\ref{T:lift}). 
One can verify that $\gamma$ is \orcon if, and only if,
$\gamma(\Sigma_{\oc}(\AAA)) \subseteq \Sigma_{\oc}(\AAA')$.
The \orpre and \lorpre homomorphisms are
defined similarly.

\begin{definition} \label{D:oplop}
Let $\gamma \colon (\GG(\AAA), \Sigma(\AAA)) \longrightarrow
(\GG(\AAA'), \Sigma(\AAA'))$ be a scale preserving homomorphism of
dimension distribution groups so that for each $n$, there is an integer
$k_n$ such that $\gamma(\GG(A_n)) \subseteq \GG(A'_{k_n})$.  We say
that $\gamma$ is {\em \orcon\/} if,
 for each $n$, the restriction
of $\gamma$ to $\GG(A_n)$ is \orcon as a map from
$\GG(A_n)$ into $\GG(A'_{k_n})$.  Similarly, $\gamma$ is
{\em \lorcon\/} if each such restriction is \lorcon.
\end{definition}

The following theorem gives in particular 
 a sufficient condition for the regular isomorphism
of systems in which all compositions are \lorcon.

\begin{thm} \label{T:clas1}
Let ${\cal A}$ and ${\cal A'}$ be systems of 
finite dimensional nest algebras 
such that all compositions of embeddings are \lorcon.
Then the following statements are equivalent.
\begin{enumerate}
\item ${\cal A}$ and ${\cal A'}$  are isomorphic with 
\lorcon  crossover maps.
\item $Inv({\cal A})$ and $Inv({\cal A'})$ are isomorphic via a 
\lorcon isomorphism.
\end{enumerate}
\end{thm}

\begin{proof}
The usual scheme of proof that  2. implies 1.
 can be completed as follows.  With the aid
of Theorem~\ref{T:lift}, obtain a lifting of a suitable restriction
\[
 \gamma_1 \colon (\GG(A_{n_1}), \Sigma(A_{n_1})) \to
   (\GG(A'_{m_1}), \Sigma(A'_{m_1})) 
\]
to a \lorcon star extendible homomorphism $\phi_1$.
In the same way, the restriction of 
$\gamma^{-1} $ to $\GG(A'_{m_1})$ may be 
lifted to a \lorcon embedding $\psi_1 $.  Since
$\mbox{$\GG(\psi_1 \circ \phi_1)$} = \GG(\beta)$, 
where $\beta$ is a
composition of the given embeddings for $\cal A$, it folows that
$\psi_1 \circ \phi_1 $ is actually \lorcon.
  From Theorem~\ref{T:r1uni} we can replace $\psi_1 $ 
by an inner unitary conjugate to
obtain $\psi_1 \circ \phi_1 = \beta $.  Continuing in this way, we
obtain a commuting diagram for the desired isomorphism.
\end{proof}

The same proof also works in the \orcon case; the next 
theorem  classifies  \orcon systems of finite dimensional nest
algebras  up to \orcon commuting diagram
isomorphism.

\begin{thm} \label{T:clas2}
Let ${\cal A}$ and ${\cal A'}$ be systems of
 finite dimensional nest algebras with
 \orcon embeddings. Then the following statements are equivalent.
\begin{enumerate}
\item ${\cal A}$ and ${\cal A'}$  are isomorphic with 
\orcon  crossover maps.
\item $Inv({\cal A})$ and $Inv({\cal A'})$ are isomorphic via an
\orcon isomorphism.
\end{enumerate}
\end{thm}

We also have parallel theorems for the \orpre and
\lorpre cases. 
The proofs are simpler in that Theorem~\ref{T:lift} 
may be replaced by the simpler
lifting argument above for \lorpre
homomorphisms.

\begin{thm} \label{T:clas3}
Let ${\cal A}$ and ${\cal A'}$ be systems of
 finite dimensional nest algebras with \orpre embeddings.
Then the following statements are equivalent.
\begin{enumerate}
\item ${\cal A}$ and ${\cal A'}$  are isomorphic with 
\orpre  crossover maps.
\item $Inv({\cal A})$ and $Inv({\cal A'})$ are isomorphic via an
\orpre isomorphism.
\end{enumerate}
\end{thm}

\begin{thm} \label{T:clas4}
Let ${\cal A}$ and ${\cal A'}$ be systems of
 finite dimensional nest algebras such that all compositions of
  embeddings are \lorpre.
Then the following statements are equivalent.
\begin{enumerate}
\item ${\cal A}$ and ${\cal A'}$  are isomorphic with 
\lorpre  crossover maps.
\item $Inv({\cal A})$ and $Inv({\cal A'})$ are isomorphic via a 
\lorpre isomorphism.
\end{enumerate}
\end{thm}

The theorems above and also Theorem~\ref{T:occlas} 
have exact counterparts for
limits of direct sums of finite dimensional nest algebras

\begin{remark} 
It is natural to ask to what extent the hypotheses of the theorems above
can be relaxed. It can be shown for example that if
$\theta \colon  K_0(A_1) \to K_0(A_2)$ is an isomorphism preserving the 
\orpre algebraic order $S_{\op} (-)$ then there is a unique (up to
conjugacy) \lorpre lifting $\phi \colon A_1 \to A_2$.
Here,  $S_{\op} (A)$ is the set of pairs $([vv^*],[v^*v])$ in 
$\Sigma_0(A) \times \Sigma_0(A)$ arising from order preserving 
 partial isometries.)  This suggests that at the level of systems,
an isomorphism $\gamma \colon K_0(\cal A) \to K_0(\cal A')$
which maps  $S_{\op} ({\cal A_1})$ to $S_{\op} ({\cal A_2})$
may lift to an isomorphism of systems. On the other hand,
a composition of locally order preserving
embeddings is not necessarily locally order preserving
 and so the usual proof is
not available. (See also Power \cite[errata]{scp92bk}.) 
This suggests that such an
isomorphism may not be sufficient for system or algebra isomorphism. In
fact this problem is already present in the case of triangular systems
where the order $S_{\op}(-)$ agrees with the algebraic order 
$S(-)$.
One way to settle the issue  would be to show 
that even for alternation algebras the invariant  
$(K_0(-), S(-))$ is not a complete invariant.
\end{remark}

Next we give applications of Theorem~\ref{T:clas2} and
\ref{T:clas3} to the classification of
algebraic limits of \orcon and \orpre systems.

Let $\sA$ belong to $\Sys (\sF_{\oc})$ or to   
$\Sys (\sF_{\op})$ with  algebraic direct limit
$A_0 = \alglim \sA$. It follows from Theorem~\ref{T:istmp} 
that we may define
$\Inv (A_0) = \Inv (\sA)$ as an invariant for star extendible
isomorphism. The next theorem shows that, in the $\Sys (\sF_{\oc})$
case, the definition
$\Sigma_{\oc}(A_0) = \Sigma_{\oc}(\sA)$
also gives an invariant and that $(\Inv (A_0),  \Sigma_{\oc}(A_0))$
is a complete invariant.

\begin{thm} \label{T:occlas}
Let $\sA,\sA'$ be \orcon systems
of finite dimensional nest algebras and let 
  $A_0, A_0'$ be their algebraic direct limits.
Then $A_0$ and $ A_0'$ are star extendibly isomorphic if, and only
if, there is an isomorphism $\gamma \colon \Inv (A_0) \to \Inv (A_0')$
such that $\gamma (\Sigma_{\oc}(A_0)) = \Sigma_{\oc}(A_0')$
\end{thm}

\begin{proof}
In view of Theorem~\ref{T:clas2}, it will be sufficient to show that 
if $A_0, A_0'$ are star extendibly isomorphic then not only 
is the induced  isomorphism 
of $\sA,\sA'$ regular, but the commuting diagram isomorphism
may be  implemented by \orcon embeddings. However this follows
immediately from Lemma~\ref{L:autooc}
\end{proof}

And, of course, there is an analogous theorem in the \orpre context,
with $\Sigma_{\op}(A_0) = \Sigma_{\op}(\sA)$
  an invariant and  $(\Inv (A_0),  \Sigma_{\op}(A_0))$
a complete invariant.

\begin{thm} \label{T:opclas}
Let $\sA,\sA'$ be \orpre systems
of finite dimensional nest algebras and let 
  $A_0, A_0'$ be their algebraic direct limits.
Then $A_0$ and $ A_0'$ are star extendibly isomorphic if, and only
if, there is an isomorphism $\gamma \colon \Inv (A_0) \to \Inv (A_0')$
such that $\gamma (\Sigma_{\op}(A_0)) = \Sigma_{\op}(A_0')$
\end{thm}

It is a consequence of Theorem 3.6 of Donsig \cite{apd98} that the
operator algebra limits of order preserving systems of finite dimensional
nest algebras are star extendibly isomorphic if, and only if, the algebraic
limit algebras are isomorphic. Combining this with Theorem~\ref{T:occlas}
it follows that the operator algebras of 
order preserving systems may be classified by
$ (\Inv (-), \Sigma_{\oc}(-))$. We anticipate that the same is true in the
\orcon case.

We conclude with some examples.

\begin{example} 
Let $A_k = T_{2^k} \otimes M_2$ and consider the embeddings  
$\phi_k \colon A_k \longrightarrow A_{k+1}$ 
 given by $\phi(a) = a \oplus a$ 
(standard type embeddings). Then
 $A = \varinjlim(A_k, \phi_k) \cong S \otimes M_2$, where
$S = \varinjlim (T_{2^k}, \sigma_k)$ is the standard upper triangular
$2^{\infty}$ limit algebra. Let
 $X_0 = \prod_{i=1}^{\infty}\{0,1\} $ with
the product topology.  The spectrum  (groupoid)
$R(C^*(S))$ is the `tails the same'
equivalence relation on $X_0$ and  
the spectrum $R(S)$ is the reverse lexicographic
sub-relation. $R(M_2) = \{0,1\} \times \{0,1\}$, the full equivalence
relation on $\{0,1\}$.  $R(A) = R(M_2) \times R(S)$, the product
relation acting on $X = \{0,1\} \times X_0$. We have 
\begin{enumerate} \renewcommand{\theenumi}{\roman{enumi}}
\renewcommand{\labelenumi}{(\theenumi)}
\item
$K_0(A) \cong C(X_0, \bbz)$, the continuous integer valued functions
  on $X_0$.
\item
$\Sigma_0(A) = \{\alpha \in K_0(A) \mid 0 \leq \alpha(x) \leq 2,
\text{ for all } x \in X_0 \}$.
\item
$\GG(A) \cong C_c(R(S), \bbz)$, the continuous integer valued
functions on $R(S)$ with compact support.
\end{enumerate}
We may make the following interpretation:
for each $p \in R(S)$, let $K_0(p) = K_0(M_2)$;
then $\GG(A)$ can be identified with
the space of compactly supported continuous functions
on $R(S)$ such that $\phi(p) \in K_0(p)$, for all $p$.
\end{example}

\begin{example}
Let $B_0$ be the unital algebraic limit algebra
$B_0 = \alglim (T_{n_k}, \phi_k)$
 where the embeddings $\phi_k$ are order preserving. Let
\[
 A_0 = \alglim (T_{n_k} \otimes M_{m_k}, \phi_k \otimes \psi_k )
\]
where each embedding $\psi_k \colon M_{m_k} \to M_{m_{k+1}}$ is 
a  unital C*-algebra injection.
Then  $A_0 = B_0 \otimes D_0$ with $D_0$ a unital ultramatricial
algebra.
As we observed at the end of 
Section 3, the maps $\phi_k \otimes \psi_k$ and their 
compositions are \orpre.
In the degenerate case $\phi_k = id \colon T_2 \to T_2$, $\GG(A_0)$ is
naturally isomorphic to 
$T_2(\mathbb{Z}) \otimes_\mathbb{Z} K_0(D_0)$. More
generally $\GG(A_0)$ can be identified naturally in terms of
$K_0(D_0)$-valued functions on the spectrum (semigroupoid) $R(B_0)$ of
$B_0$.
\end{example}

\begin{example}
 Let $B_0$ be the unital locally finite algebra
$B_0 = \alglim (M_{n_k}, \sigma_k)$ where the embeddings $\sigma_k$
are standard embeddings (of direct sum form) 
with respect to the usual matrix unit systems.
Let $A_k \subseteq M_{n_k}$ be a nest algebra whose nest projections
are ordered in a compatible way with the diagonal matrix units
and suppose that $\sigma_k(A_k) \subseteq A_{k+1}$ for all $k$.
Since the multiplicity one decomposition of $\sigma_k$
is necessarily an ordered sum, it is clear that each 
$\sigma_k$ is \orcon as a map from $A_k$
to $ A_{k+1}$.  Plainly the algebra $A_0$ contains the triangular algebra
$T^\sigma_0 = \alglim (T_{n_k}, \sigma_k)$. In fact it is the case that
{\em every} intermediate algebra $A_0'$, with the property 
$T^\sigma_0  \subseteq  A_0'  \subseteq B_0$ is necessarily of this
form.
This can be seen by a straightforward 
application of inductivity arguments (Chapter 4 of \cite{scp92bk}) 
in the setting of algebraic direct  limits;
$A_0' =  \alglim (A_k , \sigma_k)$ where $A_k' = A_0 \cap C^*(A_k)$,
and since $A_k'$ contains $T_{n_k}$, for each $k$, $A_k'$
is a finite dimensional nest algebra.
Also it follows similarly that
if $A'$ is a closed subalgebra of a UHF C*-algebra which contains
a triangular standard embedding limit algebra, then $A'$ is 
an \orcon limit of finite dimensional nest algebras.
\end{example}

\begin{example}
 In analogy with the last example consider a refinement limit
presentation of $B_0$ together with its associated triangular
refinement limit algebra $T^\rho_0$. Then it can be shown that every
intermediate algebra is an \orcon refinement limit
of finite dimensional nest algebras. This verification depends
on the inductivity of the intermediate algebras and the fact that
a refinement type embedding
between finite dimensional nest algebras 
is \orcon. 

The most transparent intermediate algebras in this case are those 
nest subalgebras of
$B_0$ determined by a finite nest of $T^\rho_0$-invariant
projections.
More generally however the intermediate algebras are 
nest subalgebras of $B_0$ determined by  a nest of 
$T^\rho_0$-invariant projections
in the weak closure of  $B_0$ in the tracial representation.
\end{example}

\begin{example}
In a similar way one obtains that an 
intermediate operator algebra $A$ satisfying
$B \subseteq A \subseteq C^*(B)$, where $B$ is an alternation algebra,
is an \orcon limit of finite dimensional nest algebras.  In fact,
the same result is valid whenever $A$ is an intermediate algebra
between a direct limit of $T_n$'s with order preserving embedding and
its enveloping \cstar algebra.
\end{example}


\providecommand{\bysame}{\leavevmode\hbox to3em{\hrulefill}\thinspace}

\end{document}